\documentclass{amsart}
\usepackage{amssymb} 
\usepackage{color}
\usepackage[all]{xy}
\usepackage{enumerate}

 \newtheorem{theorem}{Theorem}[section]

 \newtheorem{proposition}[theorem]{Proposition}
 \newtheorem{lemma}[theorem]{Lemma}
 \newtheorem{corollary}[theorem]{Corollary}
 \newtheorem{remark}[theorem]{Remark}
 \newtheorem{definition}[theorem]{Definition}
 \newtheorem{notation}[theorem]{Notation}
 \newtheorem{example}[theorem]{Example}

 \numberwithin{equation}{section}

\def\cat{{\rm Cat}_{\infty}}

\def\overcat#1{\cat\mbox{$\scriptstyle /#1$}}
\def\duo{{\rm Duo}_{\infty}}

\def\operad{{\rm Op}_{\infty}}

\def\op#1{{\rm Op}_{\infty/#1}}

\def\laxsfmoncatp#1{\mathsf{Mon}^{\mathrm{lax}}_{#1}(\cat)}
\def\laxmoncatp#1{{\rm Mon}^{\mathrm{lax}}_{#1}(\cat)}
\def\laxsfmon{\mathsf{Mon}^{\rm lax}}
\def\laxmon{{\rm Mon}^{\rm lax}}
\def\oplaxsfmon{\mathsf{Mon}^{\rm oplax}}
\def\oplaxmon{{\rm Mon}^{\rm oplax}}

\def\oplaxsfmoncatp#1{\mathsf{Mon}_{#1}^{\rm oplax}(\cat)}

\def\oplaxmoncatp#1{{\rm Mon}_{#1}^{\rm oplax}(\cat)}

\def\cart#1{{\rm Cart}_{/#1}}
\def\cocart#1{{\rm coCart}_{/#1}}
\def\oplaxcart#1{{\rm Cart}^{\rm oplax}_{/#1}}
\def\loplaxcart#1{{\rm Cart}^{\rm oplax, L}_{/#1}}
\def\laxcocart#1{{\rm coCart}^{\rm lax}_{/#1}}
\def\rlaxcocart#1{{\rm coCart}^{\rm lax, R}_{/#1}}

\begin{document}

\title
%[Duoidal $\infty$-categories, III]
%{Duoidal $\infty$-categories, III:\\
{On higher monoidal $\infty$-categories}
\author{Takeshi Torii}
\address{Department of Mathematics, 
%Faculty of Sience, 
Okayama University,
Okayama 700--8530, Japan}
\email{torii@math.okayama-u.ac.jp}
%\thanks{}

%\subjclass[2000]{Primary ; Secondary }
%\keywords{}

%\subjclass[2000]{18D10 (primary), 18D35, 18D50, 55U40 (secondary)}

\date{October 30, 2021 (version~1.1)}

\begin{abstract}
In this paper we introduce a notion of
$\mathbf{O}$-monoidal $\infty$-categories
for a finite sequence $\mathbf{O}^{\otimes}$ of $\infty$-operads,
which is a generalization of the notion of
higher monoidal categories
in the setting of $\infty$-categories.
\if0
We also introduce lax $\mathbf{O}$-monoidal functors
between coCartesian $\mathbf{O}$-monoidal $\infty$-categories
and oplax $\mathbf{O}$-monoidal functors
between Cartesian $\mathbf{O}$-monoidal $\infty$-categories.
\fi
We show that the $\infty$-category of
coCartesian $\mathbf{O}$-monoidal $\infty$-categories and
right adjoint lax $\mathbf{O}$-monoidal
functors is equivalent to the opposite
of the $\infty$-category of
Cartesian $\mathbf{O}_{\rm rev}$-monoidal $\infty$-categories
and left adjoint oplax $\mathbf{O}_{\rm rev}$-monoidal
functors,
where $\mathbf{O}^{\otimes}_{\rm rev}$
is a sequence obtained by reversing the order
of $\mathbf{O}^{\otimes}$.
\end{abstract}

\maketitle

\section{Introduction}
\label{section:introduction}

The purpose of this paper is twofold:
Firstly, we generalize the notion 
of duoidal $\infty$-categories and 
introduce higher monoidal $\infty$-categories.
Secondly, we prove that 
the $\infty$-category of 
higher monoidal $\infty$-categories 
and right adjoint lax monoidal functors
is equivalent to the opposite
of the $\infty$-category of 
higher monoidal $\infty$-categories 
and left adjoint oplax monoidal
functors.

A duoidal category is a category equipped with
two monoidal structures in which 
one is (op)lax monoidal with respect to the other.
The notion of duoidal category was
introduced in \cite{Aguiar-Mahajan} and
developed further by 
\cite{Batanin-Markl,BCZ,Bohm-Lack,Booker-Street, Street}.
In \cite{Torii1} we have introduced duoidal
$\infty$-categories which are analogues of duoidal 
categories in the setting of $\infty$-categories.
In this paper 
we will consider generalizations of duoidal $\infty$-categories
to higher monoidal $\infty$-categories.
%We will define $\infty$-categories equipped with
%two products in the sense of category of operators \cite{Barwick},
%in which one is lax with respect to the other.
%We also consider iteration 

Higher monoidal categories was also introduced
in \cite[Chapter~7]{Aguiar-Mahajan}.
Roughly speaking,
$n$-monoidal category is a category equipped
with $n$ ordered monoidal structures
related by interchange laws.
We will define a notion of
coCartesian $\mathbf{O}$-monoidal
$\infty$-categories for a finite sequence
$\mathbf{O}^{\otimes}$ of $\infty$-operads
over perfect operator categories.
This is a generalization of higher monoidal
categories and it is appropriate
setting to consider lax $\mathbf{O}$-monoidal functors.
In order to consider oplax monoidal functors,
we also introduce a notion
of Cartesian $\mathbf{O}$-monoidal $\infty$-categories.
We will show that the notion of
coCartesian $\mathbf{O}$-monoidal $\infty$-categories
is equivalent to that of
Cartesian $\mathbf{O}_{\rm rev}$-monoidal $\infty$-categories
(Corollary~\ref{cor:equivalence-O-lax-oplax}),
where $\mathbf{O}^{\otimes}_{\rm rev}$
is a sequence obtained by reversing the order
of $\mathbf{O}^{\otimes}$.
Furthermore,
by combining these notions,
we also introduce 
mixed $(\mathbf{O},\mathbf{P})$-monoidal $\infty$-categories
to discuss bilax $(\mathbf{O},\mathbf{P})$-monoidal
functors.

In \cite[Proposition~5]{Torii2}
we have shown that 
the left adjoint of a lax monoidal functor
between monoidal $\infty$-categories
is an oplax monoidal functor.
Furthermore,
Haugseng~\cite{Haugseng2}
and Hebestreit-Linskens-Nuiten~\cite{HLN} 
have independently shown that
the $\infty$-category of
$\mathcal{O}$-monoidal $\infty$-categories
and left adjoint oplax monoidal functors
is equivalent 
to the opposite of
the $\infty$-category of $\mathcal{O}$-monoidal
$\infty$-categories and right adjoint lax
monoidal functors for any
$\infty$-operad $\mathcal{O}^{\otimes}$.

Our main theorem 
is a generalization of this equivalence
to higher monoidal $\infty$-categories.
We denote by 
$\mathsf{Mon}_{\mathbf{O}}^{\rm lax,R}(\cat)$
the $\infty$-category of coCartesian $\mathbf{O}$-monoidal
$\infty$-categories and right adjoint lax
$\mathbf{O}$-monoidal functors
and
by
$\mathsf{Mon}_{\mathbf{O}_{\rm rev}}^{\rm oplax,L}(\cat)$
the $\infty$-category of Cartesian $\mathbf{O}_{\rm rev}$-monoidal
$\infty$-categories and left adjoint
oplax $\mathbf{O}_{\rm rev}$-monoidal functors.
%where $\mathbf{O}^{\otimes}_{\rm rev}$
%is a sequence of $\infty$-operads
%obtained by reversing the order of
%$\mathbf{O}^{\otimes}$.

\begin{theorem}[{Theorem~\ref{them:adjoint-monoidal-functors}}]
There exists an equivalence
\[ \mathsf{Mon}_{\mathbf{O}}^{\rm lax,R}(\cat)
   \simeq
   \mathsf{Mon}_{\mathbf{O}_{\rm rev}}^{\rm oplax,L}(\cat)^{\rm op} \]
of $\infty$-categories.
%which is given on objects by taking
%dual fibrations and on morphisms by taking adjoints fiberwise. 
\end{theorem}

The organization of this paper is as follows:
In \S\ref{section:iterated-cocartesian-fibrations}
we introduce
notions of $\mathbf{S}$-(co)Cartesian fibrations
and (op)lax $\mathbf{S}$-morphisms
for a finite sequence $\mathbf{S}$ of $\infty$-categories.
We also define mixed $(\mathbf{S},\mathbf{T})$-fibrations
and bilax $(\mathbf{S},\mathbf{T})$-morphisms
by combining $\mathbf{S}$-coCartesian fibrations
and $\mathbf{T}$-Cartesian fibrations.
In \S\ref{section:higher-monoidal-infty-categories}
we define higher monoidal $\infty$-categories
and study (op)lax monoidal functors between them.
We introduce notions of coCartesian $\mathbf{O}$-monoidal
$\infty$-categories and lax $\mathbf{O}$-monoidal
functors,
where $\mathbf{O}^{\otimes}$ is 
a finite sequence of $\infty$-operads over perfect operator categories.
We also introduce Cartesian $\mathbf{O}$-monoidal
$\infty$-categories and oplax $\mathbf{O}$-monoidal
functors.
By combining these notions,
we also define mixed $(\mathbf{O},\mathbf{P})$-monoidal
$\infty$-categories and bilax $(\mathbf{O},\mathbf{P})$-monoidal
functors.
In \S\ref{section:adjoints-monoidal-functors}
we study adjoints of (op)lax monoidal functors
between higher monoidal $\infty$-categories,
%For this purpose,
%we study adjoints of (op)lax morphisms of mixed fibrations.
%Using this,
%we prove the main theorem
and prove the main theorem
(Theorem~\ref{them:adjoint-monoidal-functors}).
\if0
which says that 
the $\infty$-category of $\mathbf{O}$-monoidal
$\infty$-categories and left adjoint oplax 
$\mathbf{O}$-monoidal functors
is equivalent to the opposite of the $\infty$-category
of $\mathbf{O}_{\rm rev}$-monoidal
$\infty$-categories and right adjoint
lax $\mathbf{O}_{\rm rev}$-monoidal functors.
\fi

%\input{notation}
%\section{notation}
%\begin{notation}\rm

\bigskip

\noindent
{\bf Notation}.
We denote by
$\cat$ the $\infty$-category of $\infty$-categories.
For an $\infty$-category $C$,
we denote by $C^{\simeq}$ the underlying $\infty$-groupoid.
We write 
${\rm Map}_C(x,y)$ for the mapping space in $C$.
For $\infty$-categories $C$ and $D$,
we write ${\rm Fun}(C,D)$
the $\infty$-categories of functors
from $C$ to $D$.

For $\infty$-categories $S$ and $T$,
we let $\pi_S: S\times T\to S$ be the projection.
For a functor $p: X\to S\times T$
of $\infty$-categories and an object $s\in S$,
we set $p_S=\pi_S\circ p$.
For an object $s\in S$,
we denote by $X_s$ the fiber of $p_S$ at $s$
and $p_s: X_s\to T$ the restriction
of $p$ to $X_s$.

Let 
$\mathbf{S}=(S_1,S_2,\ldots,S_n)$
be a finite sequence of $\infty$-categories.
We write $l(\mathbf{S})$ for
the length of $\mathbf{S}$.
For an integer $i$,
we set
$\mathbf{S}_{\ge i}=(S_i,S_{i+1},\ldots,S_n)$ and 
$\mathbf{S}_{\neq i}=(S_1,\ldots,S_{i-1},S_{i+1},\ldots,S_n)$,
etc.
We also set
$\mathbf{S}^{\rm op}=(S_1^{\rm op},S_2^{\rm op},\ldots,S_n^{\rm op})$
and
$\mathbf{S}_{\rm rev}=(S_n,S_{n-1},\ldots, S_1)$.
We denote by
$\prod \mathbf{S}$
the product $S_1\times S_2\times\cdots \times S_n$.
For finite sequences $\mathbf{S}=(S_1,\ldots,S_n)$
and $\mathbf{T}=(T_1,\ldots,T_m)$ of $\infty$-categories,
we set $[\mathbf{S},\mathbf{T}]=(S_1,S_2,\ldots,S_n,T_1,\ldots,T_m)$.

%\end{notation}

%\input{acknowledgement}
\bigskip

\noindent
{\bf Acknowledgements}.
The author would like to thank
Jonathan Beardsley
for letting him know the references
\cite{Haugseng2} and \cite{HLN}.
The author was partially supported by
by JSPS KAKENHI Grant Numbers JP17K05253.

%%%\input{monoidalcategories}
%%%\input{oplax-monoidal-categories}
%%%\input{lax-monoidal-functors}
%%%%%\input{mixed-fibration-old}
%%%\input{mixed-fibration}
%%%\input{marked-mixed-fibration}
%%%%%\input{cocart-cart-case}
%%%\input{duoidal-categories}
%%%%%\input{bilax-functors}
%%%\input{bialgebras}
%%%\input{monoidal-bicategories}
%%%\input{looping}
%%%\input{XA}
%%%\input{functoriality}
%%%\input{bicategory-to-duoidal}
%%%\input{MD}
%%%\input{proof-bicategory-to-duoidal}
%%%%%\input{cat-non-symmetric}

%\newpage
%\input{fibrations}
%\newpage

\section{Iterated (co)Cartesian fibrations
and mixed fibrations}
\label{section:iterated-cocartesian-fibrations}

In this section we introduce
notions of iterated (co)Cartesian fibrations
and mixed fibrations in a subcategory $\mathcal{Z}$ of
$\overcat{Z}$ and study their basic properties.
In \S\ref{subsection:(co)Cartesian-fibrations}
we recall the notions of
Cartesian fibrations and coCartesian fibrations
in $\mathcal{Z}$ studied in \cite{Torii1}.  
In \S\ref{subsection:mixed-fibrations}
we generalize mixed fibrations introduced
in \cite[\S3.2]{Torii1} to mixed
fibrations in $\mathcal{Z}$.
In \S\ref{subsection:iterated-(co)Cartesian-fibrations}
we define $\mathbf{S}$-coCartesian fibrations
and lax $\mathbf{S}$-morphisms,
and $\mathbf{T}$-Cartesian fibrations and 
oplax $\mathbf{T}$-morphisms
for finite sequences $\mathbf{S}$ and $\mathbf{T}$
of $\infty$-categories
by iterating the constructions 
in \S\ref{subsection:(co)Cartesian-fibrations}.
We also define mixed $(\mathbf{S},\mathbf{T})$-fibrations
by combining $\mathbf{S}$-coCartesian fibrations
and $\mathbf{T}$-Cartesian fibrations.

\subsection{Cartesian fibrations and coCartesian fibrations
in $\mathcal{Z}$}
\label{subsection:(co)Cartesian-fibrations}

In this subsection we recall the notions of
Cartesian fibrations and coCartesian fibrations
in a subcategory
of $\overcat{Z}$ studied in \cite{Torii1}. 

Let $S$ and $Z$ be $\infty$-categories.
Suppose that $\mathcal{Z}$ is a replete subcategory
of $\overcat{Z}$.
%such that
%the inclusion $\mathcal{Z}\hookrightarrow \overcat{Z}$
%is an isofibration.
In \cite[Definition~3.1]{Torii1} 
we introduced an $\infty$-category
$\mathsf{Fun}(S,\mathcal{Z})$.
In this paper we call an object of 
$\mathsf{Fun}(S,\mathcal{Z})$
a coCartesian fibration over $S$ in $\mathcal{Z}$.
We write
${\rm coCart}_{/S}(\mathcal{Z})$
for the $\infty$-category of coCartesian fibrations
over $S$ in $\mathcal{Z}$
instead of $\mathsf{Fun}(S,\mathcal{Z})$.
For the convenience of readers,
we will explicitly describe objects and morphisms
of $\mathsf{Fun}(S,\mathcal{Z})$.

%We consider a functor
%$p: X\to S\times Z$ of $\infty$-categories.
%Let $\pi_S: S\times Z\to S$ be the projection.
%For $s\in S$, we denote by
%$X_s$ the fiber of the composite $\pi_S\circ p:
%X\to S$ at $s$.
%We denote by $p_s: X_s\to Z$ 
%the restriction of $p$ to $X_s$.

\begin{definition}\rm
A coCartesian fibration over $S$ in $\mathcal{Z}$
is a functor $p: X\to S\times Z$ of $\infty$-categories
that satisfies the following conditions:

\begin{enumerate}

\item
The composite $p_S: X\to S$
is a coCartesian fibration
and the functor $p$ takes
$p_S$-coCartesian morphisms
to $\pi_S$-coCartesian morphisms.

\item
For each $s\in S$,
the restriction
$p_s: X_s\to Z$ is an object of $\mathcal{Z}$.

\item
For each morphism $s\to s'$ in $S$,
the induced functor $X_s\to X_{s'}$ over $Z$
is a morphism in $\mathcal{Z}$.

\end{enumerate}

A morphism between coCartesian fibrations 
$p: X\to S\times Z$ and $q: Y\to S\times Z$ over $S$
in $\mathcal{Z}$ is a functor 
$f: X\to Y$ over $S\times Z$ that satisfies
the following conditions:

\begin{enumerate}

\item
The functor $f$ takes $p_S$-coCartesian
morphisms to $q_S$-coCartesian morphisms.

\item
For each $s\in S$,
the restriction $f_s: X_s\to Y_s$ over $Z$
is a morphism in $\mathcal{Z}$.
\end{enumerate}

We define
\[ {\rm coCart}_{/S}(\mathcal{Z}) \]
to be the subcategory of $\overcat{S\times Z}$
with coCartesian fibrations over $S$ in $\mathcal{Z}$
and morphisms between them.
\end{definition}

By \cite[Proposition~3.2]{Torii1},
we have the following lemma.

\begin{lemma}[{\cite[Proposition~3.2]{Torii1}}]
\label{lemma:Z-(un)straightening}
There is an equivalence 
\[ {\rm Fun}(S,\mathcal{Z})\simeq
   {\rm coCart}_{/S}(\mathcal{Z}) \]
of $\infty$-categories.
\end{lemma}

We also define lax morphisms
between coCartesian fibrations over $S$ in $\mathcal{Z}$.

\begin{definition}\rm
A lax morphism between coCartesian fibrations 
$p: X\to S\times Z$ and $q: Y\to S\times Z$ over $S$
in $\mathcal{Z}$ is a functor 
$f: X\to Y$ over $S\times Z$ that satisfies
the following condition:
For each $s\in S$,
the restriction $f_s: X_s\to Y_s$ over $Z$
is a morphism in $\mathcal{Z}$.

We define
\[ {\rm coCart}_{/S}^{\rm lax}(\mathcal{Z}) \]
to be the subcategory of
$\overcat{S\times Z}$ with 
coCartesian fibrations over $S$ in $\mathcal{Z}$
and lax morphisms between them.
%Note that ${\rm coCart}_{/S}(\mathcal{Z})$
%is a wide subcategory of
%${\rm coCart}_{/S}^{\rm lax}(\mathcal{Z})$.
\end{definition}

Since an equivalence in ${\rm coCart}_{/S}^{\rm lax}(\mathcal{Z})$
is a morphism of coCartesian fibrations over $S$ in $\mathcal{Z}$,
we obtain the following lemma.

\begin{lemma}
\label{lemma:strict-lax-tilde-equivalence}
The inclusion functor
${\rm coCart}_{/S}^{\rm lax}(\mathcal{Z}) \to
   {\rm coCart}_{/S}(\mathcal{Z})$
induces an equivalence
\[ {\rm coCart}_{/S}^{\rm lax}(\mathcal{Z})^{\simeq}
   \stackrel{\simeq}{\longrightarrow}
   {\rm coCart}_{/S}(\mathcal{Z})^{\simeq}\]
of the underlying $\infty$-groupoids.
\end{lemma}
 
Dually, we call 
an object of $\mathsf{Fun}'(S,\mathcal{Z})$
in \cite[Remark~3.3]{Torii1}
a Cartesian fibration over $S$ in $\mathcal{Z}$,
and write ${\rm Cart}_{/S}(\mathcal{Z})$
for the $\infty$-category of Cartesian fibrations
over $S$ in $\mathcal{Z}$ instead of 
$\mathsf{Fun}'(S,\mathcal{Z})$.

\begin{definition}\rm
A Cartesian fibration over $S$ in $\mathcal{Z}$
is a functor $p: X\to S\times Z$ of $\infty$-categories
that satisfies the following conditions:

\begin{enumerate}

\item
The composite $p_S: X\to S$
is a Cartesian fibration and
the functor $p$ takes $p_S$-Cartesian morphisms
to $\pi_S$-Cartesian morphisms.

\item
For each $s\in S$,
the restriction
$p_s: X_s\to Z$ is an object of $\mathcal{Z}$.

\item
For each morphism $s'\to s$ in $S$,
the induced functor $X_s\to X_{s'}$ over $Z$
is a morphism in $\mathcal{Z}$.

\end{enumerate}

A morphism between Cartesian fibrations 
$p: X\to S\times Z$ and $q: Y\to S\times Z$ over $S$
in $\mathcal{Z}$ is a functor 
$f: X\to Y$ over $S\times Z$ that satisfies
the following conditions:

\begin{enumerate}

\item
The functor $f$ takes $(\pi_S\circ p)$-Cartesian
morphisms to $(\pi_S\circ q)$-Cartesian morphisms.

\item
For each $s\in S$,
the induced map $f_s: X_s\to Y_s$ over $Z$
is a morphism in $\mathcal{Z}$.
\end{enumerate}

We define 
\[ {\rm Cart}_{/S}(\mathcal{Z}) \]
to be the subcategory of $\overcat{S\times Z}$
with Cartesian fibrations over $S$ in $\mathcal{Z}$
and morphisms between them.
\end{definition}

\begin{definition}\rm
An oplax morphism between 
Cartesian fibrations 
$p: X\to S\times Z$ and $q: Y\to S\times Z$
over $S$ in $\mathcal{Z}$
is a functor $f: X\to Y$ over $S\times Z$
that satisfies the following condition:
For each $s\in S$,
the restriction $f_s: X_s\to Y_s$ over $Z$
is a morphism in $\mathcal{Z}$.  

We define
\[ {\rm Cart}_{/S}^{\rm oplax}(\mathcal{Z}) \]
to be the subcategory of
$\overcat{S\times Z}$ with
Cartesian fibrations over $S$ in $\mathcal{Z}$
and oplax morphisms between them.
%Note that ${\rm Cart}_{/S}(\mathcal{Z})$
%is a wide subcategory of
%${\rm Cart}_{/S}^{\rm oplax}(\mathcal{Z})$.
\end{definition}

In the same way as the case of coCartesian fibrations,
we obtain the following two lemmas.

\begin{lemma}[{\cite[Remark~3.3]{Torii1}}]
\label{lemma:straightening-unstraightening-Z}
There is an equivalence 
\[ {\rm Fun}(S^{\rm op},\mathcal{Z})\simeq
   {\rm Cart}_{/S}(\mathcal{Z}) \]
of $\infty$-categories.
\end{lemma}

\begin{lemma}
\label{lemma:strict-oplax-tilde-equivalence}
The inclusion functor
${\rm Cart}_{/S}^{\rm oplax}(\mathcal{Z}) \to
   {\rm Cart}_{/S}(\mathcal{Z})$
induces an equivalence
\[ {\rm Cart}_{/S}^{\rm oplax}(\mathcal{Z})^{\simeq}
   \stackrel{\simeq}{\longrightarrow}
   {\rm Cart}_{/S}(\mathcal{Z})^{\simeq}\]
of the underlying $\infty$-groupoids.
\end{lemma}

\begin{remark}
\label{remark:replete-iteration}
\rm
Notice that
$\laxcocart{S}(\mathcal{Z})$
and $\oplaxcart{S}(\mathcal{Z})$
are replete subcategory of $\overcat{S\times Z}$.
Thus, we can iterate the constructions
$\laxcocart{S}(-)$ and $\oplaxcart{S}(-)$.
\end{remark}

We need the following lemma 
in \S\ref{subsection:perfect-operator-categor}.

\begin{lemma}
\label{lemma:finite-product-creation-iteration}
If the inclusion functor
$\mathcal{Z}\hookrightarrow \overcat{Z}$ create
finite products,
then the inclusion functors
$\laxcocart{S}(\mathcal{Z})\hookrightarrow
\overcat{S\times Z}$ and 
$\oplaxcart{S}(\mathcal{Z})\hookrightarrow
\overcat{S^{\rm op}\times Z}$
also create finite products.
\end{lemma}

\proof
We will show that $\laxcocart{S}(\mathcal{Z})\hookrightarrow
\overcat{S\times Z}$ creates finite products.
The other case can be proved similarly.

It is easy to see that 
a final object of $\overcat{S\times Z}$
is also a final object of $\laxcocart{S}(\mathcal{Z})$. 
Let $p: X\to S\times Z$ and $q: Y\to S\times Z$
be objects of $\laxcocart{S}(\mathcal{Z})$.
By Lemma~\ref{lemma:Z-(un)straightening},
we regard $p$ and $q$ as objects
of ${\rm Fun}(S,\mathcal{Z})$.
By the assumption on $\mathcal{Z}$,
finite products of ${\rm Fun}(S,\mathcal{Z})$
are created in ${\rm Fun}(S,\overcat{Z})$.
Thus, 
$X\times_{S\times Z}Y\to S\times Z$ is a product of $p$ and $q$
in $\cocart{S}(\mathcal{Z})$.
We can easily verify that
$X\times_{S\times Z}Y\to S\times Z$ is a product
of $p$ and $q$ in $\laxcocart{S}(\mathcal{Z})$.
\qed

%\newpage
%\input{mixed}

%\newpage

\subsection{Mixed fibrations in $\mathcal{Z}$}
\label{subsection:mixed-fibrations}

In \cite[\S3.2]{Torii1}
we introduced mixed fibrations.
In this subsection 
we generalize them and introduce
mixed fibrations in $\mathcal{Z}$.

%Let $\mathcal{Z}$ be a subcategory of $\overcat{Z}$
%that is closed under equivalences. 
%Let $S$ and $T$ be $\infty$-categories.
%We consider a functor 
%$p: X\to S\times T\times Z$ of functors.
%Let $\pi_S: S\times T\times Z\to S$,
%$\pi_T: S\times T\times Z\to T$,
%and $\pi_{S\times T}: S\times T\times Z\to S\times T$
%be the projections.
%For $s\in S$ and $t\in T$, we denote by
%$X_s$, $X_t$, $X_{s,t}$ 
%the fibers of the composites $\pi_S\circ p: X\to S$,
%$\pi_T\circ p: X\to T$, and $\pi_{S\times T}\circ p: X\to S\times T$ 
%at $s$, $t$ and $(s,t)$, respectively.
%We denote by $p_s: X_s\to Z $,
%$p_t: X_t\to Z$, and
%$p_{(s,t)}: X_{(s,t)}\to Z$ 
%the restrictions of $p$ to $X_s$, $X_t$, and
%$X_{(s,t)}$, respectively.

\begin{definition}\rm
Let $S$ and $T$ be $\infty$-categories.
A mixed fibration over $(S,T)$ in $\mathcal{Z}$
is a functor $p: X\to S\times T\times Z$
of $\infty$-categories
that satisfies the following conditions:

\begin{enumerate}

\item
The composite $p_S: X\to S$ 
is a coCartesian fibration and the functor
$p$ takes $p_S$-coCartesian morphisms
to $\pi_S$-coCartesian morphisms.

\item
The composite $p_T: X\to T$
is a Cartesian fibration and the  
functor $p$ takes $p_T$-Cartesian
morphisms to $\pi_T$-Cartesian morphisms.

\item
For each $(s,t)\in S\times T$,
the restriction $p_{(s,t)}: X_{(s,t)}\to Z$
is an object of $\mathcal{Z}$.

\item
For each morphism $s\to s'$ in $S$
and each $t\in T$,
the induced functor $X_{(s,t)}\to X_{(s',t)}$ over $Z$
is a morphism in $\mathcal{Z}$.

\item 
For each $s\in S$ and each morphism
$t'\to t$ in $T$,
the induced functor
$X_{(s,t)}\to X_{(s,t')}$ over $Z$
is a morphism of $\mathcal{Z}$.

\end{enumerate}

A morphism between mixed fibrations
$p: X\to S\times T\times Z$ and $q: Y\to S\times T\times Z$
is a functor $f: X\to Y$ over $S\times T\times Z$
that satisfies the following conditions:

\begin{enumerate}

\item
The functor $f$ takes $p_S$-coCartesian morphisms
to $q_S$-coCartesian morphisms.

\item
The functor $f$ takes $p_T$-Cartesian morphisms
to $q_T$-Cartesian morphisms.

\item
For each $(s,t)\in S\times T$,
the restriction $f_{(s,t)}: X_{(s,t)}\to Y_{(s,t)}$
over $Z$ is a morphism in $\mathcal{Z}$.

\end{enumerate}

We define
\[ {\rm Mix}_{/(S,T)}(\mathcal{Z}) \]
to be the subcategory of $\overcat{S\times T\times Z}$
with mixed fibrations over $(S,T)$ in $\mathcal{Z}$
and morphisms between them.
\end{definition}

\begin{remark}\rm
If $p: X\to S\times T\times Z$
is a mixed fibration over $(S,T)$ in $\mathcal{Z}$,
then the composite $p_{S\times T}: X\to S\times T$
is a mixed fibration over $(S,T)$ in the sense
of \cite[Definition~3.15]{Torii1}. 
\end{remark}

We also define bilax morphisms between
mixed fibrations.

\begin{definition}\rm
A bilax morphism between mixed fibrations
$p: X\to S\times T\times Z$ and $q: Y\to S\times T\times Z$
over $(S,T)$ in $\mathcal{Z}$ is a
functor $f: X\to Y$ over $S\times T\times Z$ that
satisfies the following condition:
For each $(s,t)\in S\times T$,
the restriction $f_{(s,t)}: X_{(s,t)}\to Y_{(s,t)}$
over $Z$ is a morphism in $\mathcal{Z}$.

We define an $\infty$-category
\[ {\rm Mix}_{/(S,T)}^{\rm bilax}(\mathcal{Z}) \]
to be the subcategory of $\overcat{S\times T\times Z}$
with mixed fibrations over $(S,T)$ in $\mathcal{Z}$
and bilax morphisms between them.
\end{definition}

\begin{remark}\rm
When $\mathcal{Z}=\overcat{[0]}\simeq \cat$,
there is an equivalence
\[ {\rm Mix}_{/(S,T)}^{\rm bilax}(\mathcal{Z})\simeq
   {\rm Mfib}_{/(S,T)} \]
of $\infty$-categories, 
where the right hand side is 
the $\infty$-category of mixed fibrations
over $(S,T)$ defined in \cite[Definition~3.15]{Torii1}.
\end{remark}

We easily obtain the following two lemmas.

\begin{lemma}
There is an equivalence
\[ {\rm Mix}_{/(S,T)}(\mathcal{Z})\simeq
   {\rm Fun}(S\times T^{\rm op},\mathcal{Z})\]
of $\infty$-categories.
\end{lemma}

\begin{lemma}
The inclusion functor ${\rm Mix}_{/(S,T)}^{\rm bilax}(\mathcal{Z})\to
{\rm Mix}_{/(S,T)}(\mathcal{Z})$ induces an equivalence
\[ {\rm Mix}_{/(S,T)}^{\rm bilax}(\mathcal{Z})^{\simeq}
   \stackrel{\simeq}{\longrightarrow}
   {\rm Mix}_{/(S,T)}(\mathcal{Z})^{\simeq}\]
of the underlying $\infty$-groupoids.
\end{lemma}

Finally, 
we show that the $\infty$-category of
mixed fibrations can be described 
in terms of the constructions
of $\laxcocart{S}(-)$ and $\oplaxcart{T}(-)$.

\begin{theorem}
\label{thm:equivalence-bilax-coCart-Cart}
There are equivalences
\[ {\rm coCart}_{/S}^{\rm lax}
   ({\rm Cart}_{/T}^{\rm oplax}(\mathcal{Z}))
   \simeq
   {\rm Mix}_{/(S,T)}^{\rm bilax}(\mathcal{Z})
   \simeq
   {\rm Cart}_{/T}^{\rm oplax}
   ({\rm coCart}_{/S}^{\rm lax}(\mathcal{Z})) \]
of $\infty$-categories.
\end{theorem}

In order to prove Theorem~\ref{thm:equivalence-bilax-coCart-Cart},
we need the following lemmas. 
 
\begin{lemma}
\label{lemma:object-mixed-to-cocart-cart}
Let $p: X\to S\times T\times Z$ be a mixed
fibration over $(S,T)$ in $\mathcal{Z}$.
For each $s\in S$,
the restriction $p_s: X_s\to T\times Z$
is a Cartesian fibration over $T$ in $\mathcal{Z}$.
\end{lemma}

\proof
It suffices to show that  
$p_s: X_s\to T$ is a Cartesian 
fibration and $p_s: X_s\to T\times Z$
preserves Cartesian morphisms for each $s\in S$.
This follows from
\cite[Remark~3.18]{Torii1}.
\if0
Since $\pi_T\circ p: X\to T$
is a Cartesian fibration and $p: X\to S\times T\times Z$
preserves Cartesian fibrations,
it gives rise to a functor
$T^{\rm op}\to \overcat{S\times Z}$.
For each $s\in S$,
we have a functor
$\overcat{S\times Z}\to \overcat{Z}$
which is given by pullback along
the map $Z\simeq \{s\}\times Z\to S\times Z$.
Thus, we obtain a functor
$T^{\rm op}\to \overcat{Z}$,
which shows that
$\pi_T\circ p_s: X\to T$ is a Cartesian fibration
and that $p_s: X_s\to T\times Z$
preserves Cartesian morphisms.
\fi
\qed

\begin{lemma}
\label{lemma:object-cocart-cart-to-mixed}
If $p: X\to S\times T\times Z$ is an object
of ${\rm coCart}_{/S}^{\rm lax}({\rm Cart}_{/T}^{\rm oplax}(\mathcal{Z}))$,
then $p$ is a mixed fibration over $(S,T)$ in $\mathcal{Z}$.
\end{lemma}

\proof
It suffices to show that
$p_T: X\to T$
is a Cartesian fibration and
$p$ preserves Cartesian morphisms.
We can prove this in the same way
as the proof of \cite[Proposition~3.25]{Torii1}.
\qed

\proof[Proof of Theorem~\ref{thm:equivalence-bilax-coCart-Cart}]
We shall show that 
${\rm Mix}_{/(S,T)}^{\rm bilax}(\mathcal{Z})$ is equivalent to
${\rm coCart}_{/S}^{\rm lax}({\rm Cart}_{/T}^{\rm oplax}(\mathcal{Z}))$.
By the symmetry of the definition of
mixed fibrations,
we can prove the other equivalence similarly.

We can identify the objects of
${\rm Mix}_{/(S,T)}^{\rm bilax}(\mathcal{Z})$
with those of
${\rm coCart}_{/S}^{\rm lax}(\rm Cart)_{/T}^{\rm oplax}(\mathcal{Z})$
by Lemmas~\ref{lemma:object-mixed-to-cocart-cart} 
and \ref{lemma:object-cocart-cart-to-mixed}.
Thus, it suffices to show that 
we can also identify the morphisms
in ${\rm Mix}_{/(S,T)}^{\rm bilax}(\mathcal{Z})$
with those in
${\rm coCart}_{/S}^{\rm lax}(\rm Cart)_{/T}^{\rm oplax}(\mathcal{Z})$.

Let $f: X\to Y$ be a functor over $S\times T\times Z$,
where $p: X\to S\times T\times Z$ and $q: Y\to S\times T\times Z$
are mixed fibrations over $(S,T)$ in $\mathcal{Z}$.
The functor $f$ is a morphism of
mixed fibrations if and only if 
$f_{(s,t)}: X_{(s,t)}\to Y_{(s,t)}$ over $Z$
is a morphism in $\mathcal{Z}$ for each $(s,t)\in S\times T$.
On the other hand,
$f$ is a morphism in 
${\rm coCart}_{/S}^{\rm lax}(\rm Cart)_{/T}^{\rm oplax}(\mathcal{Z})$
if and only if 
$f_s: X_s\to Y_s$ over $T\times Z$
is a morphism in ${\rm Cart}_{/T}^{\rm oplax}(\mathcal{Z})$
for each $s\in S$.
This is equivalent to the condition
on $f$ being a morphism 
of mixed fibrations.
This completes the proof.
\qed

%\newpage
%\input{iterated}
%\newpage

\subsection{Iterated (co)Cartesian fibrations in $\mathcal{Z}$}
\label{subsection:iterated-(co)Cartesian-fibrations}

%Let 
%$\mathbf{S}=(S_1,S_2,\ldots,S_n)$
%be a finite sequence of $\infty$-categories.
%We write $l(\mathbf{S})$ for
%the length of $\mathbf{S}$.
%For an integer $i$,
%we set
%$\mathbf{S}_{\ge i}=(S_i,S_{i+1},\ldots,S_n)$, etc.
%We set
%$\prod \mathbf{S}=S_1\times S_2\times\cdots \times S_n$.

In this subsection
we will define 
$\mathbf{S}$-coCartesian fibrations
and lax $\mathbf{S}$-morphisms,
and $\mathbf{T}$-Cartesian fibrations
and oplax $\mathbf{T}$-morphisms
for finite sequences $\mathbf{S}$ and 
$\mathbf{T}$ of $\infty$-categories
by iterating the constructions 
in \S\ref{subsection:(co)Cartesian-fibrations}.
We also define mixed $(\mathbf{S},\mathbf{T})$-fibrations
by combining $\mathbf{S}$-coCartesian fibrations
and $\mathbf{T}$-Cartesian fibrations.

\begin{definition}\rm
Let 
$\mathbf{S}=(S_1,S_2,\ldots,S_n)$
be a finite sequence of $\infty$-categories.
We define an $\infty$-category 
\[ {\rm coCart}_{/\mathbf{S}}^{\rm lax}(\mathcal{Z}) \]
%\[ {\rm coCart}_{/\mathbf{S}}^{\rm lax}(\cat) \]
by induction on $l(\mathbf{S})$ as follows:
When $l(\mathbf{S})=0$,
we set
${\rm coCart}_{/\mathbf{S}}^{\rm lax}(\mathcal{Z})=\mathcal{Z}$.
%\[ {\rm coCart}_{/\mathbf{S}}^{\rm lax}(\cat)=\cat.\]
When $l(\mathbf{S})>0$,
we define 
\[ {\rm coCart}_{/\mathbf{S}}^{\rm lax}(\mathcal{Z})=
   {\rm coCart}_{/S_1}^{\rm lax}
   ({\rm coCart}_{/\mathbf{S}_{\ge 2}}^{\rm lax}(\mathcal{Z})).\]
We call
${\rm coCart}_{/\mathbf{S}}^{\rm lax}(\mathcal{Z})$
the $\infty$-category of $\mathbf{S}$-coCartesian fibrations
and lax $\mathbf{S}$-morphisms in $\mathcal{Z}$.
%Note that 
%${\rm coCart}_{/\mathbf{S}}^{\rm lax}(\cat)$
%is a full subcategory of $\overcat{\prod \mathbf{S}}$.
\end{definition}

Dually, we define an $\infty$-category of
$\mathbf{T}$-Cartesian fibration
and $\mathbf{T}$-oplax morphisms in $\mathcal{Z}$.

\begin{definition}\rm
Let 
$\mathbf{T}=(T_1,T_2,\ldots,T_n)$
be a finite sequence of $\infty$-categories.
We define an $\infty$-category 
\[ {\rm Cart}_{/\mathbf{T}}^{\rm oplax}(\mathcal{Z}) \]
by induction on $l(\mathbf{T})$ as follows:
When $l(\mathbf{T})=0$,
we set
${\rm Cart}_{/\mathbf{T}}^{\rm oplax}(\mathcal{Z})=\mathcal{Z}$.
When $l(\mathbf{T})>0$,
we define 
\[ {\rm Cart}_{/\mathbf{T}}^{\rm oplax}(\mathcal{Z})=
   {\rm Cart}_{/T_1}^{\rm oplax}
   ({\rm Cart}_{/\mathbf{T}_{\ge 2}}^{\rm oplax}(\mathcal{Z})).\]
We call 
${\rm Cart}_{/\mathbf{T}}^{\rm oplax}(\mathcal{Z})$
the $\infty$-category of $\mathbf{T}$-Cartesian fibrations
and oplax $\mathbf{T}$-morphisms in $\mathcal{Z}$.
\end{definition}

\begin{remark}\rm
The notions of $\mathbf{S}$-coCartesian fibrations
and $\mathbf{T}$-Cartesian fibrations
are related to Gray fibrations and Gray tensor
products.
See \cite[\S3]{Haugseng2} and \cite[\S2]{HLN}.

\end{remark}

Let $\mathbf{S}$ and $\mathbf{T}$ be finite
sequences of $\infty$-categories.
Now, we define an $\infty$-category
of mixed $(\mathbf{S},\mathbf{T})$-fibrations
and bilax $(\mathbf{S},\mathbf{T})$-morphisms.
 
\begin{definition}\rm
We define
\[ {\rm Mix}_{/(\mathbf{S},\mathbf{T})}^{\rm bilax}(\mathcal{Z}) \]
to be the $\infty$-category
${\rm coCart}_{/\mathbf{S}}^{\rm lax}
   ({\rm Cart}_{/\mathbf{T}}^{\rm oplax}(\mathcal{Z}))$
and call it 
the $\infty$-category of mixed $(\mathbf{S},\mathbf{T})$-fibrations
and bilax $(\mathbf{S},\mathbf{T})$-morphisms in $\mathcal{Z}$.
\end{definition}

\begin{remark}\rm
By Theorem~\ref{thm:equivalence-bilax-coCart-Cart},
there is an equivalence
\[ {\rm Mix}_{/(\mathbf{S},\mathbf{T})}^{\rm bilax}(\mathcal{Z}) 
   \simeq
   {\rm Cart}_{/\mathbf{T}}^{\rm oplax}
   ({\rm coCart}_{/\mathbf{S}}^{\rm lax}(\mathcal{Z}))\]
of $\infty$-categories.
\end{remark}

\if0
We will define subcategories
${\rm Mix}_{/(\mathbf{S},\mathbf{T})}^{\rm oplax}(\mathcal{Z})$
and
${\rm Mix}_{/(\mathbf{S},\mathbf{T})}^{\rm lax}(\mathcal{Z})$
of 
${\rm Mix}_{/(\mathbf{S},\mathbf{T})}^{\rm bilax}(\mathcal{Z})$ 
as follows.

\begin{definition}\rm
We define
\[ \begin{array}{rcl}
    {\rm Mix}_{/(\mathbf{S},\mathbf{T})}^{\rm oplax}(\mathcal{Z})
    &=&
    \cocart{\mathbf{S}}(\oplaxcart{\mathbf{T}}(\mathcal{Z})),\\[2mm]
    {\rm Mix}_{/(\mathbf{S},\mathbf{T})}^{\rm lax}(\mathcal{Z})
    &=&
    \laxcocart{\mathbf{S}}(\cart{\mathbf{T}}(\mathcal{Z})).\\
   \end{array}\]
\end{definition}
\fi

\if0
We can easily obtain the following lemma.

\begin{lemma}
Let $U$ be an $\infty$-category.
There are equivalences
\[ \begin{array}{rcl}
   {\rm Fun}(U,\oplaxcart{\mathbf{T}}(\mathcal{Z}))
   &\simeq &
   {\rm Mix}_{/(U,\mathbf{T})}^{\rm oplax}(\mathcal{Z}),\\[2mm]
   {\rm Fun}(U^{\rm op},\laxcocart{\mathbf{S}}(\mathcal{Z}))
   &\simeq&
   {\rm Mix}_{/(\mathbf{S},U)}^{\rm lax}(\mathcal{Z})\\  
   \end{array}\]
of $\infty$-categories.
\end{lemma}
\fi

Let $U$ be an $\infty$-category.
We describe the mapping space
${\rm Map}_{\cat}(U,{\rm Mix}_{/(\mathbf{S},\mathbf{T})}^{\rm bilax}(\mathcal{Z}))$
in terms of mixed fibrations in $\mathcal{Z}$.

\begin{theorem}
\label{theorem:equivalence-U-mixed-fibration}
There are natural equivalences
\[ {\rm Mix}_{/([U,\mathbf{S}],\mathbf{T})}^{\rm bilax}(\mathcal{Z})^{\simeq}
   \simeq
   {\rm Map}_{\cat}(U,
   {\rm Mix}_{/(\mathbf{S},\mathbf{T})}^{\rm bilax}(\mathcal{Z}))
   \simeq
   {\rm Mix}_{/(\mathbf{S},[U^{\rm op},\mathbf{T}])}^{\rm bilax}
   (\mathcal{Z})^{\simeq}\] 
of $\infty$-groupoids.
\end{theorem}

\proof
We shall show that 
${\rm Map}_{\cat}(U,
   {\rm Mix}_{/(\mathbf{S},\mathbf{T})}^{\rm bilax}(\mathcal{Z}))
\simeq
{\rm Mix}_{/([U,\mathbf{S}],\mathbf{T})}^{\rm bilax}(\mathcal{Z})^{\simeq}$.
The other equivalence can be proved
similarly.

By Lemmas~\ref{lemma:Z-(un)straightening}
and \ref{lemma:strict-lax-tilde-equivalence},
we have an equivalence
\[ {\rm Map}_{\cat}(U,{\rm Mix}_{/(\mathbf{S},\mathbf{T})}^{\rm bilax}
   (\mathcal{Z}))
   \simeq
   {\rm coCart}_{/U}^{\rm lax}
   ({\rm Mix}_{/(\mathbf{S},\mathbf{T})}^{\rm bilax}
   (\mathcal{Z}))^{\simeq}
\]
of $\infty$-groupoids.
By definition,
there is an equivalence
\[ 
{\rm coCart}_{/U}^{\rm lax}
   ({\rm Mix}_{/(\mathbf{S},\mathbf{T})}^{\rm bilax}
   (\cat))\simeq
   {\rm Mix}_{/([U,\mathbf{S}],\mathbf{T})}^{\rm bilax}(\cat) 
\]
of $\infty$-categories.
Combining these equivalences,
we obtain the desired equivalence.
\if0
\[ 
${\rm Fun}(U,
   {\rm Mix}_{/(\mathbf{S},\mathbf{T})}^{\rm bilax}(\cat))^{\simeq}
   \simeq
   {\rm Mix}_{/((U,\mathbf{S}),\mathbf{T})}^{\rm bilax}(\cat)^{\simeq}$. 
%\]
\[ \begin{array}{rcl}
   {\rm Fun}(U,{\rm Mix}_{/(\mathbf{S},\mathbf{T})}^{\rm bilax}
   (\cat))^{\simeq} &\simeq&
   {\rm coCart}_{/U}({\rm Mix}_{/(\mathbf{S},\mathbf{T})}^{\rm bilax}
   (\cat))^{\simeq}\\[2mm]
   &\simeq&
   {\rm coCart}_{/U}^{\rm lax}({\rm Mix}_{/(\mathbf{S},\mathbf{T})}^{\rm bilax}
   (\cat))^{\simeq}\\[2mm]
   &\simeq&
   {\rm Mix}_{/((U,\mathbf{S}),\mathbf{T})}^{\rm bilax}
   (\cat)^{\simeq}.\\
   \end{array}\]
\fi
\qed

%\bigskip

\if0
\begin{corollary}
There is an equivalence
\[ {\rm coCart}_{/\mathbf{S}}^{\rm lax}(\cat)^{\simeq}
   \simeq
   {\rm Cart}_{/\mathbf{S}^{\rm op}_{\rm rev}}^{\rm oplax}(\cat)^{\simeq}\]
of $\infty$-groupoids.
\end{corollary}

We can generalize this corollary as follows.
\fi

\begin{corollary}
\label{cor:cocart-cart-duality}
%Let $\mathcal{Z}$ be a subcategory
%of $\overcat{Z}$ closed under equivalences.
There is an equivalence
\[ {\rm coCart}_{/\mathbf{S}}^{\rm lax}
   (\mathcal{Z})^{\simeq}
   \simeq
   {\rm Cart}_{/\mathbf{S}^{\rm op}_{\rm rev}}^{\rm oplax}
   (\mathcal{Z})^{\simeq}\]
of the underlying $\infty$-groupoids. 
\end{corollary}

\begin{definition}\rm
By Corollary~\ref{cor:cocart-cart-duality},
we can see that 
a coCartesian $\mathbf{S}$-fibration canonically determines
a Cartesian $\mathbf{S}_{\rm rev}^{\rm op}$-fibration,
and vice versa. 
For a coCartesian or Cartesian $\mathbf{S}$-fibration
$p: X\to \prod\mathbf{S}$,
we denote by $p^{\vee}: X^{\vee}\to
\prod \mathbf{S}_{\rm rev}^{\rm op}$
the corresponding Cartesian or coCartesian
$\mathbf{S}_{\rm rev}^{\rm op}$-fibration,
and we call $p^{\vee}$
the dual fibration to $p$.
\end{definition}

%\bigskip

\if0
{\color{red} 以下の二つのLemmaは引用で済ませ、削除する予定。}

\begin{lemma}
\label{lemma:characterization-cocartesian-inductive}
Let $p: X\to \prod\mathbf{S}\times\prod\mathbf{T}\times Z$
be a mixed $(\mathbf{S},\mathbf{T})$-fibration
in $\mathcal{Z}$.
Let $(\mathbf{s}',\mathbf{t},z)
\in\prod\mathbf{S}_{\ge 2}\times\prod\mathbf{T}\times Z$
and 
let $e: x\to x'$ be a morphism of 
$X_{(\mathbf{s}',\mathbf{t},z)}$.
Then
$e$ is $p_{S_1}$-coCartesian 
if and only if $e$ is $p_{(\mathbf{s}',\mathbf{t},z)}$-coCartesian.
\end{lemma}

\proof
This follows from \cite[Lemma~3.19]{Torii1}.
\qed

%\bigskip

%\if0
Dually, we have
the following lemma.

\begin{lemma}
\label{lemma:characterization-cartesian-inductive}
Let $p: X\to \prod\mathbf{S}\times\prod\mathbf{T}\times Z$
be a mixed $(\mathbf{S},\mathbf{T})$-fibration in $\mathcal{Z}$.
Let $(\mathbf{s},\mathbf{t}',z)
\in\prod\mathbf{S}\times\prod\mathbf{T}_{\ge 2}\times Z$
and 
let $e: x\to x'$ be a morphism of 
$X_{(\mathbf{s},\mathbf{t}')}$.
Then
$e$ is $p_{T_1}$-Cartesian 
if and only if $e$ is $p_{(\mathbf{s},\mathbf{t}',z)}$-Cartesian.
\end{lemma}

\fi

%\newpage
%\input{perfect-operator}
%\newpage

\section{Higher monoidal $\infty$-categories}
\label{section:higher-monoidal-infty-categories}

In this section we introduce higher monoidal $\infty$-categories
and study (op)lax monoidal functors between them.
In \S\ref{subsection:perfect-operator-categor}
we define $\mathcal{O}$-monoidal $\infty$-categories
for an $\infty$-operad $\mathcal{O}^{\otimes}$
over a perfect operator category.
In \S\ref{section:o1-o2-monoidal}
we generalize 
the notion of duoidal $\infty$-categories
to that of $\mathbf{O}$-monoidal
$\infty$-categories,
where $\mathbf{O}^{\otimes}$ is 
a finite sequence of $\infty$-operads over perfect operator categories.
We also define and study 
(op)lax $\mathbf{O}$-monoidal functors between
them.
                
\subsection{$O$-monoidal $\infty$-categories 
over perfect operator categories}
%(op)lax $\mathcal{O}$-monoidal functors}
%for an $\infty$-operad $\mathcal{O}$
%over a perfect opertor category}
\label{subsection:perfect-operator-categor}

In this subsection we define 
$\mathcal{O}$-monoidal $\infty$-categories
for an $\infty$-operad 
$\mathcal{O}^{\otimes}$ over a  perfect operator category, 
and introduce lax and oplax $\mathcal{O}$-monoidal 
functors between $\mathcal{O}$-monoidal 
$\infty$-categories.
We also consider mixed $(\mathcal{O},\mathcal{P})$-monoidal
$\infty$-categories and mixed $(\mathcal{O},\mathcal{P})$-monoidal
functors.

First, we briefly recall the notion of
$\infty$-operads over a perfect operator category
introduced in \cite{Barwick}.
Let $\Phi$ be a prefect operator category 
in the sense of \cite[Definitions~1.2 and 4.6]{Barwick}.
%By \cite[Theorem~5.10]{Barwick},
%we have a canonical monad $T$ on $\Phi$. 
Associated to $\Phi$,
we have the Leinster category $\Lambda(\Phi)$
equipped with collections of inert morphisms and active morphisms.
%The objects of $\Lambda(\Phi)$ are those of $\Phi$.
%The morphisms for two objects $I,J$ of $\Phi$ are given by
%\[ {\rm Hom}_{\Lambda(\Phi)}(J,I)={\rm Hom}_{\Phi}(J,TI).\]
%According to \cite[Definition~7.1]{Barwick},
%we say that a morphism $J\to I$ of $\Lambda(\Phi)$
%is inert if the correponding morphism $J\to TI$
%in $\Phi$ has the property that 
%the natural morphism
%\[ J\times_{TI}I\longrightarrow I \]
%is an isomorphism.
%We denote by $\Lambda^{\dagger}(\Phi)$
%the set of inert morphisms of $\Lambda(\Phi)$.
According to \cite[Definition~7.8]{Barwick},
an $\infty$-operad over $\Phi$
is a functor
$p: \mathcal{O}^{\otimes}\to \Lambda(\Phi)$
of $\infty$-categories satisfying the following conditions:
\begin{enumerate}
\item
For every inert morphism $\phi:I\to J$ of $\Lambda(\Phi)$
and every object $x\in \mathcal{O}^{\otimes}_I$,
there is a $p$-coCartesian morphism $x\to y$
in $\mathcal{O}^{\otimes}$ covering $\phi$.
\item
For any objects $I,J\in\Lambda(\Phi)$,
any objects $x\in \mathcal{O}^{\otimes}_I$ and 
$y\in \mathcal{O}^{\otimes}_J$,
any morphism $\phi:I\to J$ of $\Lambda(\Phi)$,
and any $p$-coCartesian morphisms
$\{y\to y_j|\ j\in |J|\}$ lying over
the inert morphisms $\{\rho_j: J\to \{j\}|\ j\in |J|\}$,
the induced map 
\[ {\rm Map}_{\mathcal{O}^{\otimes}}^{\phi}(x,y)\longrightarrow
   \prod_{j\in |J|}{\rm Map}_{\mathcal{O}^{\otimes}}^{\rho_j\circ\phi}(x,y_j) \]
is an equivalence.
\item
For any object $I\in\Lambda(\Phi)$, 
the $p$-coCartesian morphisms lying over the inert
morphisms $\{I\to \{i\}|\ i\in |I|\}$ together
induce an equivalence
\[ \mathcal{O}^{\otimes}_I\longrightarrow 
   \prod_{i\in |I|}\mathcal{O}^{\otimes}_{i} .\]
\end{enumerate}

For an $\infty$-operad 
$p:\mathcal{O}^{\otimes}\to \Lambda(\Phi)$
over $\Phi$,
a morphism of $\mathcal{O}^{\otimes}$
is said to be inert if it is a $p$-coCartesian morphism
over an inert morphism,
and active if it covers an active morphism
of $\Lambda(\Phi)$.
A morphism of $\infty$-operads over $\Phi$
between $\mathcal{O}^{\otimes}\to\Lambda(\Phi)$
and $\mathcal{P}^{\otimes}\to\Lambda(\Phi)$
is a functor $f: \mathcal{O}^{\otimes}\to \mathcal{P}^{\otimes}$
over $\Lambda(\Phi)$ that preserves inert morphisms.

\begin{example}
\rm
Let $\mathbf{F}$ be the category of finite sets.
By \cite[Example~4.9.3]{Barwick},
$\mathbf{F}$ is a perfect operator category.
By \cite[Example~6.5]{Barwick},
$\Lambda(\mathbf{F})$ 
is equivalent to the category ${\rm Fin}_*$ 
of pointed finite sets.
By \cite[Example~7.9]{Barwick},
the notion of $\infty$-operads over $\mathbf{F}$
coincides with that of Lurie's 
$\infty$-operads
in \cite[Chapter~2]{Lurie2}. 
\end{example}

\begin{example}
\rm
Let $\mathbf{O}$ be the category of ordered finite sets.
By \cite[Example~4.9.2]{Barwick},
$\mathbf{O}$ is a perfect operator category.
By \cite[Example~6.6]{Barwick},
$\Lambda(\mathbf{O})$ is equivalent to $\Delta^{\rm op}$.
The notion of $\infty$-operads over $\mathbf{O}$
coincides with that of non-symmetric $\infty$-operads
in \cite[\S3]{Gepner-Haugseng}. 
\end{example}

\if0
\begin{definition}\rm
Let $p: X^{\otimes}\to \mathcal{O}_1^{\otimes}\times
(\mathcal{O}_2^{\otimes})^{\rm op}$ be a categorical
fibration.
We say that $p: X^{\otimes}\to \mathcal{O}_1^{\otimes}\times
(\mathcal{O}_2^{\otimes})^{\rm op}$ is 
an $\mathcal{O}_1$-$\mathcal{O}_2$-duoidal
$\infty$-category if
it is a virtual duoidal $\mathcal{O}_1$-$\mathcal{O}_2$-duoidal
$\infty$-category and 
it satisfies the following conditions;
\begin{enumerate}
\item
For each $I_1\in\Lambda(\Phi_1)$,
the restriction $p_{I_1,\bullet}:
X^{\otimes}_{I_1,\bullet}\to (\mathcal{O}_2^{\otimes})^{\rm op}$
is a Cartesian fibration.
\item
For each $I_2\in\Lambda(\Phi_2)$<
the restriction $p_{\bullet,I_2}:
X^{\otimes}_{\bullet,I_2}\to \mathcal{O}_1^{\otimes}$
is a coCartesian fibration.
\end{enumerate}
\end{definition}
\fi

For a perfect operator category $\Phi$,
we denote by $\operad^{\Phi}$  
the $\infty$-category of $\infty$-operads
over $\Phi$.
By definition,
a perfect operator category $\Phi$
has a final object $*$
(see \cite[Definition~1.2]{Barwick}
for the definition of operator categories).
For an $\infty$-operad $\mathcal{O}^{\otimes}\to \Lambda(\Phi)$,
we denote by $\mathcal{O}$
the fiber $\mathcal{O}^{\otimes}_*$ at $*\in \Lambda(\Phi)$
and say that it is the $\infty$-category of 
colors of $\mathcal{O}^{\otimes}$. 

We define a coCartesian 
$\mathcal{O}$-monoidal $\infty$-category
for an $\infty$-operad $\mathcal{O}^{\otimes}$ over $\Phi$.

\begin{definition}\rm
Let 
$\mathcal{C}^{\otimes}\to\mathcal{O}^{\otimes}$
be a morphism of $\infty$-operads over $\Phi$.
We say that $\mathcal{C}$ is a coCartesian 
$\mathcal{O}$-monoidal $\infty$-category if 
the functor $\mathcal{C}^{\otimes}\to\mathcal{O}^{\otimes}$
is a coCartesian fibration.
When it is clear from the context,
we simply say that it is an $\mathcal{O}$-monoidal
$\infty$-category. 
\end{definition}

Let $\mathcal{C}^{\otimes}\to\mathcal{O}^{\otimes}$
be a coCartesian $\mathcal{O}$-monoidal $\infty$-category.
%We say that $\mathcal{C}$ is the underlying 
%$\infty$-category.
%For $T\in\mathcal{O}^{\otimes}$,
%we denote by $\mathcal{C}^{\otimes}_T$
%the fiber at $T\in\mathcal{O}^{\otimes}$.
For each color 
$x\in\mathcal{O}$,
we write $\mathcal{C}_x$ for $\mathcal{C}^{\otimes}_x$
and call it the underlying $\infty$-category of
$\mathcal{C}^{\otimes}$ over $x$.

By the same proof of
\cite[Proposition~2.1.2.12]{Lurie2},
we have the following lemma.

\begin{lemma}
\label{lemma:cocartesian-operads}
Let $p:\mathcal{O}^{\otimes}\to \Lambda(\Phi)$
be an $\infty$-operad over $\Phi$,
and let $f: \mathcal{C}^{\otimes}\to\mathcal{O}^{\otimes}$
be a coCartesian fibration.
Then the following conditions are equivalent\,\mbox{\rm :}
\begin{enumerate} 
\item[{\rm (a)}]
The composite $q: \mathcal{C}^{\otimes}\stackrel{f}{\to}
\mathcal{O}^{\otimes}\stackrel{p}{\to}\Lambda(\Phi)$
is an $\infty$-operad over $\Phi$.
\item[{\rm (b)}]
For every $x\in \mathcal{O}^{\otimes}_I$,
the $p$-coCartesian morphisms $\{x\to x_i|\ i\in |I|\}$
covering the inert morphisms $\{I\to \{i\}|\ i\in |I|\}$
induce an equivalence
$\mathcal{C}^{\otimes}_x\to \prod_{i\in |I|}
\mathcal{C}_{x_i}$
of $\infty$-categories.
\end{enumerate}
\end{lemma}

A coCartesian $\mathcal{O}$-monoidal 
$\infty$-category $\mathcal{C}$ 
is equipped with multiplications indexed by
active morphisms in $\mathcal{O}^{\otimes}$
as in \cite[Remark~2.1.2.16]{Lurie2}.
Let $x,y\in \mathcal{O}^{\otimes}$
with $p(x)=I$ and $p(y)=1$.
We have a unique active morphism $a: I\to 1$  
in $\Lambda(\Phi)$.
Let $\{x\to x_i|\ i\in |I|\}$ be
inert morphisms of $\mathcal{O}^{\otimes}$
covering the inert morphisms $\{I\to \{i\}|\ i\in |I|\}$.
For every $\theta\in {\rm Map}_{\mathcal{O}^{\otimes}}^a(x,y)$,
we have a multiplication map
\[ \otimes_{\theta}:
   \prod_{i\in |I|}\mathcal{C}_{x_i}\stackrel{\simeq}{\longleftarrow}
   \mathcal{C}^{\otimes}_{x}\stackrel{\theta_*}{\longrightarrow}
   \mathcal{C}_{y},\]
where $\theta_*$ is a functor
induced by $\theta$ by using the coCartesian fibration
$\mathcal{C}^{\otimes}\to\mathcal{O}^{\otimes}$.

We denote by $\op{\mathcal{O}^{\otimes}}^{\Phi}$
the $\infty$-category 
of $\infty$-operads over $\mathcal{O}^{\otimes}$.

\begin{definition}\rm
%Let $\mathcal{C}$ and $\mathcal{D}$
%be $\mathcal{O}$-monoidal $\infty$-categories.
A lax $\mathcal{O}$-monoidal functor
between coCartesian $\mathcal{O}$-monoidal $\infty$-categories
is a morphism in $\op{\mathcal{O}^{\otimes}}^{\Phi}$. 
Furthermore, if %a morphism in $\op{\mathcal{O}^{\otimes}}^{\Phi}$ 
%between $\mathcal{O}$-monidal $\infty$-categories
it preserves coCartesian morphisms,
then we say that it is a (strong) $\mathcal{O}$-monoidal functor.
\end{definition}

\if0
For an $\infty$-category $S$,
we write 
\[ (\cocart{S})_{\rm lax},\quad
   (\cart{S})_{\rm oplax} \]
for the full subcategories of 
$\overcat{S}$ 
spanned by coCartesian fibrations and Cartesian fibrations,
respectively.
We also write
\[ \cocart{S},\quad
   \cart{S} \]
for the wide subcategories of
$(\cocart{S})_{\rm lax}$ and
$(\cart{S})_{\rm oplax}$
with morphisms that preserve coCartesian morphisms
and Cartesian morphisms, 
respectively.
\fi

\begin{definition}\rm
We define an $\infty$-category
\[ \laxsfmoncatp{\mathcal{O}} \] 
to be the subcategory of 
$\laxcocart{\mathcal{O}^{\otimes}}(\cat)$
with coCartesian $\mathcal{O}$-monoidal $\infty$-categories 
as objects
and lax $\mathcal{O}$-monoidal functors as morphisms.
\if0
the full subcategory 
of $\op{\mathcal{O}^{\otimes}}^{\Phi}$
spanned by $\mathcal{O}$-monoidal $\infty$-categories.
\fi
%We call it the $\infty$-category of $\mathcal{O}$-monoidal
%$\infty$-categories and lax $\mathcal{O}$-monoidal
%functors.

We write
\[ \mathsf{Mon}_{\mathcal{O}}(\cat) \]
for the wide subcategory of $\laxsfmoncatp{\mathcal{O}}$
with (strong) $\mathcal{O}$-monoidal functors
as morphisms.
\end{definition}
%We say that
%the fiber $\mathcal{C}=\mathcal{C}^{\otimes}_{\langle 1\rangle}$
%over $\langle 1\rangle\in\fin$
%is the underlying $\infty$-category 
%for an $\mathcal{O}$-monoidal $\infty$-category
%$\mathcal{C}^{\otimes}$.

\if0
We note that there is a functor
\[ U: \laxmoncatp{\mathcal{O}}\longrightarrow \cat,\]
which associates the underlying $\infty$-category
to an $\mathcal{O}$-monoidal $\infty$-category. 
For simplicity,
we say that 
$\mathcal{C}$ is an $\mathcal{O}$-monoidal
$\infty$-category 
when $\mathcal{C}^{\otimes}\to\mathcal{O}^{\otimes}$
is an $\mathcal{O}$-monoidal $\infty$-category,
and that 
$F: \mathcal{C}\to \mathcal{D}$
is a lax $\mathcal{O}$-monoidal functor
when a map $F: \mathcal{C}^{\otimes}\to\mathcal{D}^{\otimes}$
over $\mathcal{O}^{\otimes}$
is a lax $\mathcal{O}$-monoidal functor.
\fi

%\begin{definition}\rm
Let $\mathcal{X}$ be an $\infty$-category with finite products.
Let $p: \mathcal{O}^{\otimes}\to\Lambda(\Phi)$
be an $\infty$-operad over a perfect operator
category $\Phi$.
We say that a functor $F: \mathcal{O}^{\otimes}\to \mathcal{X}$
is an $\mathcal{O}$-monoid object of $\mathcal{X}$
if the functor
\[ F(x)\longrightarrow \prod_{i\in |I|} F(x_i) \]
is an equivalence in $\mathcal{X}$
for any $x\in\mathcal{O}^{\otimes}$,
where $p(x)=I$ and 
$\{x\to x_i|\ i\in |I|\}$ are $p$-coCartesian morphisms
lying over inert morphisms
$\{I\to \{i\}|\ i\in |I|\}$.
%\end{definition}
We denote by
\[ {\rm Mon}_{\mathcal{O}}(\mathcal{X}) \]
the full subcategory of
${\rm Fun}(\mathcal{O}^{\otimes},\mathcal{X})$
spanned by $\mathcal{O}$-monoid objects.

\if0
\begin{definition}\rm
Let $p: \mathcal{O}^{\otimes}\to\Lambda(\Phi)$
be an $\infty$-operad over a perfect operator
category $\Phi$.
A functor $C: \mathcal{O}^{\otimes}\to \cat$
is said to be an $\mathcal{O}$-monoid object
if for any $x\in\mathcal{O}^{\otimes}$ with $p(x)=I$,
the map
\[ C(x)\longrightarrow \prod_{i\in |I|} C(x_i) \]
is an equivalence,
where $\{x\to x_i|\ i\in |I|\}$
is a set of $p$-coCartesian morphisms
lying over the set of inert morphisms
$\{I\to \{i\}|\ i\in |I|\}$.
\end{definition}
\fi

To a coCartesian $\mathcal{O}$-monoidal 
$\infty$-category $\mathcal{C}$, 
by the straightening functor 
of coCartesian fibrations
\cite[\S3.2]{Lurie1},
we can associate
a functor $C: \mathcal{O}^{\otimes}\to\cat$,
which is an $\mathcal{O}$-monoid object in $\cat$. 
Furthermore,
by Lemma~\ref{lemma:cocartesian-operads},
there is an equivalence
\[ \mathsf{Mon}_{\mathcal{O}}(\cat) \simeq
   \mathrm{Mon}_{\mathcal{O}}(\cat)\]
of $\infty$-categories.

Let $C: \mathcal{O}^{\otimes}\to \cat$
be an $\mathcal{O}$-monoid object
corresponding to a coCartesian $\mathcal{O}$-monoidal
$\infty$-category $\mathcal{C}^{\otimes}\to\mathcal{O}^{\otimes}$.
By the unstraightening functor
of Cartesian fibrations~\cite[\S3.2]{Lurie1},
we obtain a Cartesian fibration
$(\mathcal{C}^{\otimes})^{\vee}\to (\mathcal{O}^{\otimes})^{\rm op}$.
We call it a Cartesian $\mathcal{O}$-monoidal $\infty$-category.

We write $\mathcal{C}^{\vee}$
for $(\mathcal{C}^{\otimes})^{\vee}_1$.
For a color $x\in \mathcal{O}$,
we write $\mathcal{C}^{\vee}_x$ for 
$(\mathcal{C}^{\otimes})^{\vee}_x$
and call it the underlying $\infty$-category of
$(\mathcal{C}^{\otimes})^{\vee}$ over $x$.
Note that there are equivalences
\[ \mathcal{C}_x\simeq C(x)\simeq \mathcal{C}^{\vee}_x\]
for any $x\in\mathcal{O}$.

\if0
Let $R:\cat\to\cat$ be a functor
which associate to an $\infty$-category its opposite.
Since $R$ preserves finite products,
the composite
$RC: \mathcal{O}^{\otimes}\stackrel{C}{\to}\cat\stackrel{R}{\to}\cat$
is also an $\mathcal{O}$-monoid object of $\cat$.
We denote by
\[ R\mathcal{C}^{\otimes}\longrightarrow \mathcal{O}^{\otimes} \]
an $\mathcal{O}$-monoidal $\infty$-category
associated to the $\mathcal{O}$-monoid object $RC$.
Note that the underlying $\infty$-category
of $R\mathcal{C}^{\otimes}$ is equivalent
to the opposite $\mathcal{C}^{\rm op}$
of the underlying $\infty$-category $\mathcal{C}$.
%We denote by $R\mathcal{C}$
%the underlying $\infty$-category of $R\mathcal{C}^{\otimes}$.
\fi
 
We say that a morphism of $(\mathcal{C}^{\otimes})^{\vee}$
is inert if
it is a Cartesian morphism
over an inert morphism of $(\mathcal{O}^{\otimes})^{\rm op}$.

\begin{definition}\rm
Let $(\mathcal{C}^{\otimes})^{\vee}\to(\mathcal{O}^{\otimes})^{\rm op}$ 
and $(\mathcal{D}^{\otimes})^{\vee}\to(\mathcal{O}^{\otimes})^{\rm op}$
be Cartesian $\mathcal{O}$-monoidal $\infty$-categories.
An oplax $\mathcal{O}$-monoidal functor
from $\mathcal{C}^{\vee}$ to $\mathcal{D}^{\vee}$
is a functor
$(\mathcal{C}^{\otimes})^{\vee}\to (\mathcal{D}^{\otimes})^{\vee}$
over $(\mathcal{O}^{\otimes})^{\rm op}$
which preserves inert morphisms.
%Cartesian morphisms
%over inert morphisms of $(\mathcal{O}^{\otimes})^{\rm op}$.

We define an $\infty$-category 
\[ \oplaxsfmoncatp{\mathcal{O}}\]
to be the subcategory of 
$\oplaxcart{(\mathcal{O}^{\otimes})^{\rm op}}(\cat)$
%$\overcat{(\mathcal{O}^{\otimes})^{\rm op}}$
with Cartesian $\mathcal{O}$-monoidal $\infty$-categories
as objects and oplax $\mathcal{O}$-monoidal functors
as morphisms.
%We call it the $\infty$-category of 
%$\mathcal{O}$-monoidal $\infty$-categories
%and oplax $\mathcal{O}$-monoidal functors.
\end{definition}

Let 
$\mathcal{Z}$ be a replete subcategory
of $\overcat{Z}$.
We assume that the inclusion functor
$\mathcal{Z}\hookrightarrow \overcat{Z}$ 
creates finite products.
%is closed under
%finite products and equivalences in $\overcat{Z}$.
As in \cite[\S3.1]{Torii1},
we will introduce $\infty$-categories
$\laxsfmon_{\mathcal{O}}(\mathcal{Z})$
and
$\oplaxsfmon_{\mathcal{O}}(\mathcal{Z})$
by generalizing $\laxsfmon_{\mathcal{O}}(\cat)$
and $\oplaxsfmon_{\mathcal{O}}(\cat)$.

\if0
We write 
\[ {\rm Mon}_{\mathcal{O}}(\mathcal{X}) \]
for the full subcategory of ${\rm Fun}(\mathcal{O}^{\otimes},\mathcal{C})$
spanned by $\mathcal{O}$-monoid objects in $\mathcal{C}$.
\fi

\begin{definition}\rm
We define an $\infty$-category
\[ \laxsfmon_{\mathcal{O}}(\mathcal{Z}) \]
to be a subcategory of 
$\laxcocart{\mathcal{O}^{\otimes}}(\mathcal{Z})$
as follows.
An object of $\laxsfmon_{\mathcal{O}}(\mathcal{Z})$
is an object
$p: \mathcal{C}^{\otimes}\to \mathcal{O}^{\otimes}\times Z$
of $\laxcocart{\mathcal{O}^{\otimes}}(\mathcal{Z})$
that satisfies the following condition:
%\begin{enumerate}
%\item[]
For each $x\in\mathcal{O}^{\otimes}$,
the Segal morphism
\[ \mathcal{C}^{\otimes}_x\stackrel{\simeq}{\longrightarrow}
   \prod_{i\in |I|}\!\!{}^Z\ \mathcal{C}^{\otimes}_{x_i},\]
is an equivalence in $\mathcal{Z}$,
where the right hand side
is a product in $\overcat{Z}$.
%$p(x)=I$, and
%$\{x\to x_i|\ i\in |I|\}$ are
%inert morphisms of $\mathcal{O}^{\otimes}$
%over inert morphisms
%$\{I\to \{i\}|\ i\in |I|\}$.
%\end{enumerate}

A morphism of $\laxsfmon_{\mathcal{O}}(\mathcal{Z})$
between objects $p: \mathcal{C}^{\otimes}\to 
\mathcal{O}^{\otimes}\times Z$
and $q: \mathcal{D}^{\otimes}\to\mathcal{O}^{\otimes}\times Z$
is a morphism $f: \mathcal{C}^{\otimes}\to\mathcal{D}^{\otimes}$
in $\laxcocart{\mathcal{O}^{\otimes}}(\mathcal{Z})$
%that satisfies the following condition:
%
%\begin{enumerate}
%
%\item[]
%The functor $f$ 
that takes $p_{\mathcal{O}^{\otimes}}$-coCartesian
morphisms over inert morphisms of $\mathcal{O}^{\otimes}$
to $q_{\mathcal{O}^{\otimes}}$-coCartesian morphisms.
%
%\end{enumerate}
%\end{definition}

%\begin{definition}\rm
Dually, we define an $\infty$-category
\[ \oplaxsfmon_{\mathcal{O}}(\mathcal{Z}) \]
to be the subcategory of 
$\oplaxcart{(\mathcal{O}^{\otimes})^{\rm op}}(\mathcal{Z})$
as follows.
An object of $\oplaxsfmon_{\mathcal{O}}(\mathcal{Z})$
is an object 
$p: \mathcal{C}^{\otimes}\to (\mathcal{O}^{\otimes})^{\rm op}\times Z$
of $\oplaxcart{(\mathcal{O}^{\otimes})^{\rm op}}(\mathcal{Z})$
that satisfies the following condition:
%
%\begin{enumerate}
%
%\item[]
For each $y\in(\mathcal{O}^{\otimes})^{\rm op}$,
the Segal morphism
\[ \mathcal{C}^{\otimes}_y\stackrel{\simeq}{\longrightarrow}
   \prod_{i\in |I|}\!\!{}^Z\ \mathcal{C}^{\otimes}_{y_i}\]
is an equivalence in $\mathcal{Z}$,
where the right hand side
is a product in $\overcat{Z}$.
%$p(y)=I$, and
%$\{y\to y_i|\ i\in |I|\}$ are
%inert mrophsisms of $\mathcal{O}^{\otimes}$
%over inert morphisms
%$\{I\to \{i\}|\ i\in |I|\}$.
%
%\end{enumerate}

A morphism of $\oplaxsfmon_{\mathcal{O}}(\mathcal{Z})$
between objects 
$p: \mathcal{C}^{\otimes}\to (\mathcal{O}^{\otimes})^{\rm op}\times Z$
and $q: \mathcal{D}^{\otimes}\to(\mathcal{O}^{\otimes})^{\rm op}\times Z$
is a morphism
$f: \mathcal{C}^{\otimes}\to\mathcal{D}^{\otimes}$
in $\oplaxcart{(\mathcal{O}^{\otimes})^{\rm op}}(\mathcal{Z})$
%that satisfies the following condition:
%
%\begin{enumerate}
%
%\item[]
%The functor $f$ carries 
that takes
$p_{(\mathcal{O}^{\otimes})^{\rm op}}$-Cartesian
morphisms over inert morphisms of $(\mathcal{O}^{\otimes})^{\rm op}$
to $q_{(\mathcal{O}^{\otimes})^{\rm op}}$-Cartesian morphisms.
%
%\end{enumerate}
\end{definition}

We will show that the underlying $\infty$-groupoids
of $\laxsfmon_{\mathcal{O}}(\mathcal{Z})$
and $\oplaxsfmon_{\mathcal{O}}(\mathcal{Z})$
are equivalent to that of ${\rm Mon}_{\mathcal{O}}(\mathcal{Z})$.

\begin{lemma}
\label{lemma:lax-oplax-equivalence-groupoids}
The inclusion functors
induce equivalences
\[ \laxsfmon_{\mathcal{O}}(\mathcal{Z})^{\simeq}
   \stackrel{\simeq}{\longrightarrow} 
   {\rm Mon}_{\mathcal{O}}(\mathcal{Z})^{\simeq}
   \stackrel{\simeq}{\longleftarrow}
   \oplaxsfmon_{\mathcal{O}}(\mathcal{Z})^{\simeq}\]
of the underlying $\infty$-groupoids.
\end{lemma}

\proof
By Lemma~\ref{lemma:strict-lax-tilde-equivalence},
the inclusion functor
$\laxcocart{\mathcal{O}^{\otimes}}(\mathcal{Z})
\to \cocart{\mathcal{O}^{\otimes}}(\mathcal{Z})$
induces an equivalence
$\laxcocart{\mathcal{O}^{\otimes}}(\mathcal{Z})^{\simeq}
\stackrel{\simeq}{\to} \cocart{\mathcal{O}^{\otimes}}
(\mathcal{Z})^{\simeq}$
of the underlying $\infty$-groupoids.
This induces an equivalence
between $(\laxsfmon_{\mathcal{O}}(\mathcal{Z}))^{\simeq}$
and $({\rm Mon}_{\mathcal{O}}(\mathcal{Z}))^{\simeq}$.
The other equivalence can be proved similarly
by using Lemma~\ref{lemma:strict-oplax-tilde-equivalence}.
\qed

\bigskip

In \cite[Definition~4.18]{Torii1}
we introduced the $\infty$-category
${\rm Duo}_{\infty}^{\rm bilax}$
of duoidal $\infty$-categories and bilax monoidal
functors.
By \cite[Theorem~4.20 and Remark~4.21 ]{Torii1},
we showed that there are equivalences
\[ \laxsfmon_{\Delta^{\rm op}}(\oplaxsfmon_{\Delta^{\rm op}}(\cat))\simeq 
   {\rm Duo}_{\infty}^{\rm bilax}\simeq
   \oplaxsfmon_{\Delta^{\rm op}}(\laxsfmon_{\Delta^{\rm op}}(\cat)) \]
of $\infty$-categories.
In the following of this subsection
we will generalize these equivalences.

%Let $\mathcal{Z}$ be a subcategory
%of $\overcat{Z}$
%that is closed under finite products and equivalences,
%where $Z$ is an $\infty$-category.

%\begin{remark}\rm
%If $\mathcal{Z}$ is a replete subcategory of $\overcat{Z}$
%in which the inclusion functor
%$\mathcal{Z}\hookrightarrow\overcat{Z}$
%creates finite products, 
%is closed under finite products and equivalences,
%then 

First, 
we prove the following lemma,
which guarantees that we can iterate the constructions
$\mathsf{Mon}_{\mathcal{O}}^{\rm lax}(-)$
and 
$\mathsf{Mon}_{\mathcal{O}}^{\rm oplax}(-)$.

\begin{lemma}
\if0
Let $\mathcal{Z}$ be a replete subcategory of $\overcat{Z}$
in which the inclusion functor
$\mathcal{Z}\hookrightarrow\overcat{Z}$
creates finite products. 
\fi
The $\infty$-category
$\laxsfmon_{\mathcal{O}}(\mathcal{Z})$
is a replete subcategory of $\overcat{\mathcal{O}^{\otimes}\times Z}$,
and
the inclusion functor
$\laxsfmon_{\mathcal{O}}(\mathcal{Z})
\hookrightarrow\overcat{\mathcal{O}^{\otimes}\times Z}$
creates finite products.
Similarly,
$\oplaxsfmon_{\mathcal{O}}(\mathcal{Z})$
is a replete subcategory of 
$\overcat{(\mathcal{O}^{\otimes})^{\rm op}\times Z}$,
and the inclusion functor
$\oplaxsfmon_{\mathcal{O}}(\mathcal{Z})
\hookrightarrow\overcat{(\mathcal{O}^{\otimes})^{\rm op}\times Z}$
creates finite products.
%respectively.
%\end{remark}
\end{lemma}

\proof
By Remark~\ref{remark:replete-iteration},
we can easily see that
$\laxsfmon_{\mathcal{O}}(\mathcal{Z})$
and 
$\oplaxsfmon_{\mathcal{O}}(\mathcal{Z})$
are replete subcategories.
Using Lemma~\ref{lemma:finite-product-creation-iteration},
we can verify that 
the inclusion functors
$\laxsfmon_{\mathcal{O}}(\mathcal{Z})
\hookrightarrow\overcat{\mathcal{O}^{\otimes}\times Z}$
and
$\oplaxsfmon_{\mathcal{O}}(\mathcal{Z})
\hookrightarrow\overcat{(\mathcal{O}^{\otimes})^{\rm op}\times Z}$
create finite products.
\qed

\bigskip

\if0
\proof
It is easy to verify that
$\mathsf{Mon}_{\mathcal{O}}^{\rm lax}(\mathcal{Z})$
and
$\mathsf{Mon}_{\mathcal{O}^{\rm op}}^{\rm oplax}(\mathcal{Z})$
are replete subcategories.
We will show that if
\qed
\fi

%\bigskip

Next,
we define mixed $(\mathcal{O},\mathcal{P})$-monoidal
$\infty$-categories in $\mathcal{Z}$ and bilax 
$(\mathcal{O},\mathcal{P})$-monoidal functors
between them.

\begin{definition}
\label{definition:mixed-monoidal-infty-category}
\rm
Let $\mathcal{O}^{\otimes}\to \Lambda(\Phi)$ 
and $\mathcal{P}^{\otimes}\to \Lambda(\Psi)$
be $\infty$-operads over perfect
operator categories over
$\Phi$ and $\Psi$, respectively.
A mixed $(\mathcal{O},\mathcal{P})$-monoidal $\infty$-category 
in $\mathcal{Z}$ is a mixed fibration
$p: \mathcal{C}^{\otimes}\to \mathcal{O}^{\otimes}\times
(\mathcal{P}^{\otimes})^{\rm op}\times Z$
over $(\mathcal{O}^{\otimes},(\mathcal{P}^{\otimes})^{\rm op})$
in $\mathcal{Z}$ that satisfies the following conditions:

\begin{enumerate}

\item
For each $x\in \mathcal{O}^{\otimes}$,
the Segal morphism
\[ \mathcal{C}^{\otimes}_x
   \longrightarrow
   \prod_{i\in |I|}\!\!{}^{(\mathcal{P}^{\otimes})^{\rm op}\times Z}\
   \mathcal{C}^{\otimes}_{x_i} \]
is an equivalence
in $\oplaxcart{(\mathcal{P}^{\otimes})^{\rm op}}(\mathcal{Z})$,
where the right hand side
is a product in 
the $\infty$-category $\overcat{(\mathcal{P}^{\otimes})^{\rm op}
\times Z}$.

\item
For each $y\in (\mathcal{P}^{\otimes})^{\rm op}$,
the Segal morphism
\[ \mathcal{C}^{\otimes}_y
   \longrightarrow
   \prod_{j\in |J|}\!\!{}^{\mathcal{O}^{\otimes}\times Z}\
   \mathcal{C}^{\otimes}_{y_j} \]
is an equivalence 
in $\laxcocart{\mathcal{O}^{\otimes}}(\mathcal{Z})$,
where the right hand side 
is a product in $\overcat{\mathcal{O}^{\otimes}\times Z}$.

\end{enumerate}

A bilax $(\mathcal{O},\mathcal{P})$-monoidal functor
between mixed $(\mathcal{O},\mathcal{P})$-monoidal
$\infty$-categories $p: \mathcal{C}^{\otimes}\to 
\mathcal{O}^{\otimes}\times (\mathcal{P}^{\otimes})^{\rm op}
\times Z$ and $q: \mathcal{D}^{\otimes}\to \mathcal{O}^{\otimes}
\times (\mathcal{P}^{\otimes})^{\rm op}\times Z$ in $\mathcal{Z}$
is a bilax 
$(\mathcal{O}^{\otimes},(\mathcal{P}^{\otimes})^{\rm op})$-morphism $f$ 
%of mixed fibrations
%over $(\mathcal{O}^{\otimes},(\mathcal{P}^{\otimes})^{\rm op})$ 
%in $\mathcal{Z}$
that satisfies the following conditions:

\begin{enumerate}

\item
The functor $f$ takes 
$p_{\mathcal{O}^{\otimes}}$-coCartesian
morphisms over inert morphisms of $\mathcal{O}^{\otimes}$
to $q_{\mathcal{O}^{\otimes}}$-coCartesian morphisms.

\item
The functor $f$ takes
$p_{(\mathcal{P}^{\otimes})^{\rm op}}$-Cartesian
morphisms over inert morphisms of $(\mathcal{P}^{\otimes})^{\rm op}$
to $q_{(\mathcal{P}^{\otimes})^{\rm op}}$-Cartesian
morphisms.

\end{enumerate}

We define
\[ \mathsf{Mon}_{(\mathcal{O},\mathcal{P})}^{\rm bilax}(\mathcal{Z}) \]
to be the subcategory of
${\rm Mix}_{(\mathcal{O}^{\otimes},(\mathcal{P}^{\otimes})^{\rm op})}^{\rm bilax}
(\mathcal{Z})$ with
mixed $(\mathcal{O},\mathcal{P})$-monoidal $\infty$-categories
and 
bilax $(\mathcal{O},\mathcal{P})$-monoidal functors.
\end{definition}

\begin{theorem}
\label{theorem:lax-oplax-mon-equivalence}
There are equivalences
\[ \laxsfmon_{\mathcal{O}}(\oplaxsfmon_{\mathcal{P}}(\mathcal{Z}))
   \simeq
   \mathsf{Mon}_{(\mathcal{O},\mathcal{P})}^{\rm bilax}(\mathcal{Z})
   \simeq
   \oplaxsfmon_{\mathcal{P}}(\laxsfmon_{\mathcal{O}}(\mathcal{Z})) \]
of $\infty$-categories.
\end{theorem}

\proof
We will prove 
$\mathsf{Mon}_{(\mathcal{O},\mathcal{P})}^{\rm bilax}(\mathcal{Z})
   \simeq
\laxsfmon_{\mathcal{O}}(\oplaxsfmon_{\mathcal{P}}(\mathcal{Z}))$.
The other equivalence can be proved similarly.

By Theorem~\ref{thm:equivalence-bilax-coCart-Cart},
we have an equivalence
\[    {\rm Mix}_{/(\mathcal{O}^{\otimes},(\mathcal{P}^{\otimes})^{\rm op})}
    ^{\rm bilax}(\mathcal{Z})
   \simeq    
   {\rm coCart}_{/\mathcal{O}^{\otimes}}^{\rm lax}
   ({\rm Cart}_{/(\mathcal{P}^{\otimes})^{\rm op}}^{\rm oplax}(\mathcal{Z})).\]
Under this equivalence,
we can easily verify that
an $(\mathcal{O},\mathcal{P})$-monoidal
$\infty$-category in $\mathcal{Z}$ corresponds to
an object of 
$\mathsf{Mon}_{\mathcal{O}}^{\rm lax}
(\mathsf{Mon}_{\mathcal{P}}^{\rm oplax}(\mathcal{Z}))$,
and vice versa,
by using the same argument 
in the proof of \cite[Proposition~4.13]{Torii1}.
Furthermore, 
we see that
a bilax $(\mathcal{O},\mathcal{P})$-monoidal
functor between $(\mathcal{O},\mathcal{P})$-monoidal
$\infty$-categories corresponds to 
a morphism of $\mathsf{Mon}_{\mathcal{O}}^{\rm lax}
(\mathsf{Mon}_{\mathcal{P}}^{\rm oplax}(\mathcal{Z}))$
by using the argument in the 
proof of \cite[Theorem~4.20]{Torii1}.
\qed

\subsection
%{$\mathcal{O}_1$-$\mathcal{O}_2$-monoidal $\infty$-categories}
%{$\mathcal{O}_1$-$\mathcal{O}_2$-monoidal $\infty$-categories}
{Higher monoidal $\infty$-categories}
\label{section:o1-o2-monoidal}

In this subsection we generalize 
the notion of duoidal $\infty$-categories
to that of $\mathbf{O}$-monoidal
$\infty$-categories,
where $\mathbf{O}^{\otimes}$ is 
a finite sequence of $\infty$-operads over perfect operator categories.
%in the sense of \cite{Barwick}.

\if0
Let $\mathbf{S}=(S_1,S_2,\ldots,S_n)$
be a finite sequence of $\infty$-categories.
We write $l(\mathbf{S})$ for the length
of the sequence $\mathbf{S}$.
For $1\le i\le l(\mathbf{S})$,
we write $\mathbf{S}_{\ge i}$
for the sequence 
$(S_i,S_{i+1},\ldots,S_n)$, etc.
\fi

Let $\mathbf{O}^{\otimes}=(\mathcal{O}_1^{\otimes},%\mathcal{O}_2^{\otimes},
\ldots,\mathcal{O}_n^{\otimes})$
be a finite sequence of $\infty$-operads over
perfect operator categories
$\Phi_1,%\Phi_2,
\ldots,\Phi_n$, respectively.
We set $\mathbf{O}=(\mathcal{O}_1,\ldots,\mathcal{O}_n)$.

First, we define an $\infty$-category
of coCartesian $\mathbf{O}$-monoidal $\infty$-categories and 
lax $\mathbf{O}$-monoidal functors 
by induction on $l(\mathbf{O})$.

\begin{definition}\rm
When $l(\mathbf{O})=0$,
we set 
$\laxsfmon_{\mathbf{O}}(\cat)=\cat$.
%\[ \laxsfmon_{\mathbf{O}}(\cat)=\laxsfmon_{\mathcal{O}_1}(\cat). \]
When $l(\mathbf{O})>0$,
we define
\[ \laxsfmon_{\mathbf{O}}(\cat)=
   \laxsfmon_{\mathcal{O}_1}(\laxsfmon_{\mathbf{O}_{\ge 2}}(\cat)). \]
%where $\mathbf{O}_{\ge 2}=
%(\mathcal{O}_2^{\otimes},\mathcal{O}_3^{\otimes},
%\ldots,\mathcal{O}_n^{\otimes})$.
%We call it the $\infty$-category of 
%coCartesian $\mathbf{O}$-monoidal
%$\infty$-categories and lax $\mathbf{O}$-monoidal functors.
\end{definition}

For the convenience of readers, 
we explicitly describe objects
and morphisms of $\laxsfmon_{\mathbf{O}}(\cat)$.  

\if0
For a finite sequence $\mathbf{S}=(S_1,\ldots,S_n)$
of $\infty$-categories, we set
\[ \prod \mathbf{S}= S_1
   \times S_2\times\cdots
   \times S_n.\]
We denote by
\[ \pi_i: \prod\mathbf{S}\to S_i \]
the projection for $1\le i\le l(\mathbf{O})$.
\fi

\if0
For a categorical fibration
$p: \mathcal{C}^{\otimes}\to \prod\mathbf{O}^{\otimes}$,
we set
$p_i=\pi_i\circ p: \mathcal{C}^{\otimes}\to \mathcal{O}_i^{\otimes}$
for $1\le i\le l(\mathbf{O})$.
\fi

\begin{definition}\rm
Let $p: \mathcal{C}^{\otimes}\to \prod\mathbf{O}^{\otimes}$
be a functor of $\infty$-categories.

When $l(\mathbf{O})=0$,
we say that any functor $p: \mathcal{C}^{\otimes}\to [0]$ 
is a coCartesian $\mathbf{O}$-monoidal
$\infty$-category, 
where $[0]$ is a final object of $\cat$.
\if0
if it satisfies the following conditions:
\begin{enumerate}
\item
The map $p$ is a coCartesian fibration.
\item
For each $x\in \mathcal{O}_1^{\otimes}$,
the Segal morphism
\[ \mathcal{C}^{\otimes}_x\stackrel{\simeq}{\longrightarrow}
   \prod_{i\in |I|} \mathcal{C}^{\otimes}_{x_i} \]
is an equivalence of $\infty$-categories. 
\end{enumerate}
\fi
A lax $\mathbf{O}$-monoidal functor
between coCartesian $\mathbf{O}$-monoidal $\infty$-category
is a morphism in $\overcat{[0]}\simeq\cat$.

When $l(\mathbf{O})>0$,
we say that $p$ is a coCartesian $\mathbf{O}$-monoidal
$\infty$-category 
if it satisfies the following conditions:

\begin{enumerate}

\item 
The composite map 
$p_{\mathcal{O}_1^{\otimes}}: 
\mathcal{C}^{\otimes}\to \mathcal{O}_1^{\otimes}$ 
is a coCartesian fibration,
and $p$ takes $p_{\mathcal{O}_1^{\otimes}}$-coCartesian morphisms
to $\pi_{\mathcal{O}_1^{\otimes}}$-coCartesian morphisms.

\item
For each $x\in \mathcal{O}_1^{\otimes}$,
the restriction $\mathcal{C}^{\otimes}_x\to
\prod\mathbf{O}_{i\ge 2}^{\otimes}$
is a coCartesian $\mathbf{O}_{i\ge 2}$-monoidal $\infty$-category.

\item
For each morphism
$x\to x'$ in $\mathcal{O}_1^{\otimes}$,
the induced functor
$\mathcal{C}^{\otimes}_x\to \mathcal{C}^{\otimes}_{x'}$
over $\prod\mathbf{O}^{\otimes}_{\ge 2}$
is a lax $\mathbf{O}_{\ge 2}$-monoidal
functor.

\item
For each $x\in\mathcal{O}^{\otimes}_1$,
the Segal morphism
\[ \mathcal{C}^{\otimes}_x\stackrel{\simeq}{\longrightarrow}
%   \prod_{i\in |I|}
   \prod_{i\in |I|}\!\!{}^{\prod\mathbf{O}^{\otimes}_{\ge 2}}\ 
    \mathcal{C}^{\otimes}_{x_i},\]
is an equivalence 
of $\mathbf{O}_{\ge 2}$-monoidal $\infty$-categories,
%in $\overcat{\prod\mathbf{O}^{\otimes}_{\ge 2}}$,
where the right hand side
is a product in $\overcat{\prod\mathbf{O}^{\otimes}_{\ge 2}}$.
%$p(x)=I$, and
%$\{x\to x_i|\ i\in |I|\}$ are
%inert morphisms of $\mathcal{O}^{\otimes}_1$
%over inert morphisms
%$\{I\to \{i\}|\ i\in |I|\}$.

\end{enumerate} 
A lax $\mathbf{O}$-monoidal functor
between coCartesian $\mathbf{O}$-monoidal
$\infty$-categories $p: \mathcal{C}^{\otimes}\to \prod\mathbf{O}^{\otimes}$
and $q:\mathcal{D}^{\otimes}\to\prod\mathbf{O}^{\otimes}$
%in coCartesian fibration form
is a functor $h: \mathcal{C}^{\otimes}\to \mathcal{D}^{\otimes}$
over $\prod\mathbf{O}^{\otimes}$ that satisfies
the following conditions:

\begin{enumerate}

\item
The functor $h$ takes $p_{\mathcal{O}_1^{\otimes}}$-coCartesian 
morphisms over
inert morphisms of $\mathcal{O}_1^{\otimes}$ to 
$q_{\mathcal{O}_1^{\otimes}}$-coCartesian morphisms.

\item
For each $x\in \mathcal{O}_1^{\otimes}$,
the induced functor $h_x: \mathcal{C}^{\otimes}_x
\to \mathcal{D}^{\otimes}_x$ over
$\prod\mathbf{O}^{\otimes}_{\ge 2}$
is a lax $\mathbf{O}_{\ge 2}$-monoidal functor. 

\end{enumerate}
\end{definition}

Unwinding the definition,
we obtain the following theorem.

\begin{theorem}
\label{thm:o-lax-monoidal-cacartesian-form}
The $\infty$-category
$\laxsfmon_{\mathbf{O}}(\cat)$
is a subcategory of $\overcat{\prod\mathbf{O}^{\otimes}}$
with coCartesian $\mathbf{O}$-monoidal $\infty$-categories
as objects and
lax $\mathbf{O}$-monoidal functors as morphisms. 
\end{theorem}

\begin{example}\rm
Let $p: \mathcal{C}^{\otimes}\to\Lambda(\mathbf{F})$
be a symmetric monoidal $\infty$-category.
We obtain a coCartesian $\mathbf{O}$-monoidal $\infty$-category
for any finite sequence $\mathbf{O}$ of $\infty$-operads
over $\mathbf{F}$
by taking pullback of $p$ along
the map $\prod\mathbf{O}^{\otimes}\to \prod\Lambda(\mathbf{F})
\to \Lambda(\mathbf{F})$,
where the second map is
a smash product functor
in \cite[Notation~2.2.5.9]{Lurie2}.
\end{example}

\begin{example}\rm
Let $\mathbf{O}^{\otimes}=
(\mathcal{O}_1^{\otimes},\ldots,\mathcal{O}_n^{\otimes})$ be a finite
sequence of $\infty$-operads over $\mathbf{F}$.
We denote by $\otimes\mathbf{O}^{\otimes}$
the Boardman-Vogt tensor product
$\mathcal{O}_1^{\otimes}\otimes\cdots\otimes\mathcal{O}_n^{\otimes}$.
Let $p: \mathcal{C}^{\otimes}\to \otimes \mathbf{O}^{\otimes}$
be an $\otimes \mathbf{O}^{\otimes}$-monoidal
$\infty$-category.
We obtain a coCartesian $\mathbf{O}$-monoidal $\infty$-category
by taking pullback of $p$ along
the map $\prod\mathbf{O}^{\otimes}\to
\otimes\mathbf{O}^{\otimes}$.
\end{example}

\begin{example}\rm
Let $\mathbb{E}_k^{\otimes}$ be the little $k$-cubes operad.
Suppose that $\mathcal{C}^{\otimes}\to \mathbb{E}_{m+n}^{\otimes}$ is a 
presentable $\mathbb{E}_{m+n}$-monoidal
$\infty$-category and that
$A$ is an $\mathbb{E}_{m+n}$-algebra object of $\mathcal{C}$.
Then there exists a coCartesian
$(\mathbb{E}_m,\mathbb{E}_n)$-monoidal $\infty$-category
$\mathcal{D}^{\otimes}\to \mathbb{E}_m^{\otimes}
\times\mathbb{E}_n^{\otimes}$
such that the underlying $\infty$-category
$\mathcal{D}^{\otimes}_{(\langle 1\rangle,\langle 1\rangle)}$
is equivalent to the $\infty$-category
${\rm Mod}^{\mathbb{E}_m}_A(\mathcal{C})$
of $\mathbb{E}_m$-$A$-modules. 
See \cite{Torii3} for more details.
\end{example}

Dually,
we define an $\infty$-category 
of Cartesian $\mathbf{O}$-monoidal $\infty$-category
and oplax $\mathbf{O}$-monoidal functors.

%Next, we shall define an $\infty$-category
%of $\mathbf{O}$-monoidal $\infty$-categories and 
%oplax $\mathbf{O}$-monoidal functors 
%by induction on $l(\mathbf{O})$.

\begin{definition}\rm
When $l(\mathbf{O})=0$,
we set 
$\oplaxsfmon_{\mathbf{O}}(\cat)=\cat$.
%\[ \oplaxsfmon_{\mathbf{O}}(\cat)=\oplaxsfmon_{\mathcal{O}_1}(\cat). \]
When $l(\mathbf{O})>0$,
we define
\[ \oplaxsfmon_{\mathbf{O}}(\cat)=
   \oplaxsfmon_{\mathcal{O}_1}
   (\oplaxsfmon_{\mathbf{O}_{\ge 2}}(\cat)). \]
%where $\mathbf{O}_{\le 2}=
%(\mathcal{O}_1^{\otimes},\mathcal{O}_2^{\otimes},
%\ldots,\mathcal{O}_{n-1}^{\otimes})$.
We call it the $\infty$-category of Cartesian $\mathbf{O}$-monoidal
$\infty$-categories and oplax $\mathbf{O}$-monoidal
functors.
\end{definition}

For finite sequences $\mathbf{O}^{\otimes}$ and 
$\mathbf{P}^{\otimes}$ of $\infty$-operads
over perfect operator categories,
we also define an $\infty$-category of 
$(\mathbf{O},\mathbf{P})$-monoidal $\infty$-categories
and bilax $(\mathbf{O},\mathbf{P})$-monoidal functors.

\begin{definition}\rm
We define 
\[ \mathsf{Mon}_{(\mathbf{O},\mathbf{P})}
   ^{\mbox{\scriptsize bilax}}(\cat)=
   \laxsfmon_{\mathbf{O}}(\oplaxsfmon_{\mathbf{P}}(\cat)). \]
Note that there is an equivalence
$\mathsf{Mon}_{(\mathbf{O},\mathbf{P})}
   ^{\mbox{\scriptsize bilax}}(\cat)\simeq
   \oplaxsfmon_{\mathbf{P}}(\laxsfmon_{\mathbf{O}}(\cat))$
by Theorem~\ref{theorem:lax-oplax-mon-equivalence}.
\end{definition}

\if0
\begin{remark}\rm
We have
\[ {\rm Mon}^{\mathbf{P}\mbox{\scriptsize -lax}}
   _{\mathbf{Q}\mbox{\scriptsize -oplax}}(\cat)
   \simeq
   {\rm Mon}^{\mathbf{Q}\mbox{\scriptsize -oplax}}
   ({\rm Mon}^{\mathbf{P}\mbox{\scriptsize -lax}}(\cat)).\]
\end{remark}
\fi

\if0
\begin{definition}
Let $p: \mathcal{C}\to \prod\mathbf{O}^{\otimes}
\times\prod\mathbf{P}^{\otimes}$
be a categorical fibration.
We 
\begin{enumerate}

\item
$p_{\mathbf{O}}: \mathcal{C}^{\otimes}\to \prod\mathbf{O}^{\otimes}$
is a $\mathbf{O}$-monoidal $\infty$-category,
and $p$ is a 

\end{enumerate} 
\end{definition}
\fi

Now, we show that 
mixed higher monoidal $\infty$-categories
can canonically be identified with 
coCartesian higher monoidal $\infty$-categories or
Cartesian higher monoidal $\infty$-categories. 

\begin{theorem}
There are natural equivalences
\[  \laxsfmon_{[\mathbf{P}_{\rm rev},\mathbf{O}]}
    (\cat)^{\simeq}
    \simeq
   {\rm Mon}_{(\mathbf{O},\mathbf{P})}^{\rm bilax}
   (\cat)^{\simeq}
   \simeq
   \oplaxsfmon_{[\mathbf{O}_{\rm rev},\mathbf{P}]}
   (\cat)^{\simeq}\]
of the underlying $\infty$-groupoids.
\end{theorem}

\if0
\begin{lemma}
Let $\mathcal{Q}^{\otimes}$ be an $\infty$-operad
over a perfect operator category.
%Let $\mathbf{P}^{\otimes}$
%and $\mathbf{Q}^{\otimes}$ be 
%finite sequences of $\infty$-operads
%over perfect operator categories.
The equivalence in Theorem~\ref{theorem:equivalence-U-mixed-fibration}
restricts to natural equivalences
\[ \mathsf{Mon}_{([\mathcal{Q},\mathbf{O}],\mathbf{P})}
    ^{\rm bilax}(\cat)^{\simeq}
   \simeq 
   {\rm Mon}_{\mathcal{Q}}
   (\mathsf{Mon}_{(\mathbf{O},\mathbf{P})}^{\rm bilax}(\cat))^{\simeq}
   \simeq
   \mathsf{Mon}_{(\mathbf{O},[\mathcal{Q},\mathbf{P}])}
    ^{\rm bilax}(\cat)^{\simeq}
\]
of the underlying $\infty$-groupoids.
\end{lemma}
\fi

\proof
\if0
We will prove that 
$(\mathsf{Mon}_{([\mathcal{O},\mathbf{P}],\mathbf{Q})}
    ^{\rm bilax}(\cat))^{\simeq}
   \simeq 
   {\rm Mon}_{\mathcal{O}}
   (\mathsf{Mon}_{(\mathbf{P},\mathbf{Q})}^{\rm bilax}(\cat))^{\simeq}$.
\fi
For an $\infty$-operad $\mathcal{Q}^{\otimes}$ 
over a perfect operator category,
it suffices to show that 
the equivalence in Theorem~\ref{theorem:equivalence-U-mixed-fibration}
restricts to natural equivalences
\[ \mathsf{Mon}_{([\mathcal{Q},\mathbf{O}],\mathbf{P})}
    ^{\rm bilax}(\cat)^{\simeq}
   \simeq 
   {\rm Mon}_{\mathcal{Q}}
   (\mathsf{Mon}_{(\mathbf{O},\mathbf{P})}^{\rm bilax}(\cat))^{\simeq}
   \simeq
   \mathsf{Mon}_{(\mathbf{O},[\mathcal{Q},\mathbf{P}])}
    ^{\rm bilax}(\cat)^{\simeq}.
\]
%of the underlying $\infty$-groupoids.
By Lemma~\ref{lemma:lax-oplax-equivalence-groupoids},
we have equivalences
\[ \begin{array}{rcl}
   {\rm Mon}_{\mathcal{Q}}
   (\mathsf{Mon}_{(\mathbf{O},\mathbf{P})}^{\rm bilax}(\cat))^{\simeq}
   &\simeq&
   \mathsf{Mon}_{\mathcal{Q}}^{\rm lax}
   (\mathsf{Mon}_{(\mathbf{O},\mathbf{P})}^{\rm bilax}(\cat))^{\simeq}\\[2mm]
   &\simeq&
   \mathsf{Mon}_{([\mathcal{Q},\mathbf{O}],\mathbf{P})}^{\rm bilax}
   (\cat)^{\simeq}.\\   
   \end{array}\] 
The other equivalence can be proved similarly.
\qed

\begin{corollary}
\label{cor:equivalence-O-lax-oplax}
The equivalence in Corollary~\ref{cor:cocart-cart-duality}
restricts to a natural equivalence
\[ \laxsfmon_{\mathbf{O}}(\cat)^{\simeq}
   \simeq
   \oplaxsfmon_{\mathbf{O}_{\rm rev}}(\cat)^{\simeq} \]
of the underlying $\infty$-groupoids.
\end{corollary}

By Corollary~\ref{cor:equivalence-O-lax-oplax},
a coCartesian $\mathbf{O}$-monoidal $\infty$-category
canonically determines
a Cartesian $\mathbf{O}_{\rm rev}$-monoidal $\infty$-category,
and vice versa, 
by taking dual fibrations.

Next,
we give a fiberwise criterion
on bilax monoidal functors.

\begin{proposition}
\label{prop:criterion-(op)lax-monoidal-functor}
Let $f: \mathcal{C}^{\otimes}\to \mathcal{D}^{\otimes}$
be a bilax 
$(\mathbf{O}^{\otimes},(\mathbf{P}^{\otimes})^{\rm op})$-morphism
between $(\mathbf{O},\mathbf{P})$-monoidal 
$\infty$-categories.
\if0
$p: \mathcal{C}^{\otimes}\to 
\prod\mathbf{O}^{\otimes}\times
\prod(\mathbf{P}^{\otimes})^{\rm op}$
and
$q: \mathcal{D}^{\otimes}\to 
\prod\mathbf{O}^{\otimes}\times
\prod(\mathbf{P}^{\otimes})^{\rm op}$.
\fi
Then $f$ is a bilax $(\mathbf{O},\mathbf{P})$-monoidal
functor if and only if 
the following two conditions holds\,\mbox{\rm :}
\begin{enumerate}
\item[\mbox{\rm (a)}]
$f_{(\mathbf{s}',\mathbf{t})}$
is a lax $\mathcal{O}_i$-monoidal
functor for each $(\mathbf{s}',\mathbf{t})\in
\prod\mathbf{O}_{\neq i}^{\otimes}\times\prod
(\mathbf{P}^{\otimes})^{\rm op}$.
\item[\mbox{\rm (b)}]
$f_{(\mathbf{s},\mathbf{t}')}$
is an oplax $\mathcal{P}_j$-monoidal functor
for each $(\mathbf{s},\mathbf{t}')\in
\prod\mathbf{O}^{\otimes}\times\prod
(\mathbf{P}^{\otimes}_{\neq j})^{\rm op}$.
\end{enumerate}
\end{proposition}

\proof
It is clear that if $f$ is a bilax 
$(\mathbf{O},\mathbf{P})$-monoidal functor,
then (a) and (b) holds.
Thus, 
it suffices to show that 
if $f$ satisfies (a) and (b),
then $f$ is a bilax $(\mathbf{O},\mathbf{P})$-monoidal
functor.

We may assume that $l(\mathbf{O})+l(\mathbf{P})>1$.
We suppose that $l(\mathbf{O})>0$.
The case $l(\mathbf{P})>0$ can be proved 
similarly.

We let $p_1: \mathcal{C}^{\otimes}\to \mathcal{O}_1^{\otimes}$
and $q_1: \mathcal{D}^{\otimes}\to \mathcal{O}_1^{\otimes}$
be the projection.
Then $f$ is bilax $(\mathbf{O},\mathbf{P})$-monoidal
if (1) $f$ takes $p_1$-coCartesian morphisms
over inert morphisms of $\mathcal{O}_1^{\otimes}$
to $q_1$-coCartesian morphisms,
and (2) $f_x:\mathcal{C}^{\otimes}_x
\to\mathcal{D}^{\otimes}_x$ is a bilax 
$(\mathbf{O}_{\ge 2},\mathbf{P})$-monoidal functor
for every $x\in\mathcal{O}_1^{\otimes}$.
%By Lemmas~\ref{lemma:characterization-cocartesian-inductive},
By \cite[Lemma~3.19]{Torii1},
(1) is equivalent to
the condition that 
$f_{\mathbf{y}}$ is lax $\mathcal{O}_1$-monoidal
for every $\mathbf{y}\in \prod\mathbf{O}_{\ge 2}^{\otimes}
\times\prod(\mathbf{P}^{\otimes})^{\rm op}$.
Thus, by induction on $l(\mathbf{O})$,
$f$ is bilax $(\mathbf{O},\mathbf{P})$-monoidal
if $f$ satisfies (a) and (b).
\qed

\if0
Next, we suppose that $l(\mathbf{O})=0$
and $l(\mathbf{O})>1$.
We let $p_1: \mathcal{C}^{\otimes}\to (\mathcal{P}_1^{\otimes})^{\rm op}$
and $q_1: \mathcal{D}^{\otimes}\to (\mathcal{P}_1^{\otimes})^{\rm op}$
be the projection.
Then $f$ is lax $(\mathbf{O},\mathbf{P})$-monoidal
if (1) $f$ takes $p_1$-Cartesian morphisms
over inert morphisms of $(\mathcal{P}_1^{\otimes})^{\rm op}$
to $q_1$-Cartesian morphisms,
and (2) $f_x:\mathcal{C}^{\otimes}_x
\to\mathcal{D}^{\otimes}_x$ is a bilax 
$(\mathbf{O},\mathbf{P}_{\ge 2})$-monoidal functor
for every $x\in (\mathcal{P}_1^{\otimes})^{\rm op}$.
By Lemma~\ref{lemma:characterization-cartesian-inductive},
(1) is equivalent to
the condition that 
$f_{\mathbf{y}}$ is oplax $\mathcal{P}_1$-monoidal
for every $\mathbf{y}\in \prod\mathbf{O}^{\otimes}
\times\prod(\mathbf{P}_{\ge 2}^{\otimes})^{\rm op}$.
Thus, by induction on $l(\mathbf{P})$,
$f$ is bilax $(\mathbf{O},\mathbf{P})$-monoidal
if $\epsilon(f)$ is a morphism of $\mathcal{M}(\mathbf{O},\mathbf{P})$.
\qed
\fi

\if0
Thus, if $f$ is a bilax $(\mathbf{O},\mathbf{P})$-monoidal
functor,
then $f_{(\mathbf{x}',\mathbf{y})}$ 
is a lax $\mathcal{O}$-monoidal functor 
for any $(\mathbf{x}',\mathbf{y})\in
\prod\mathbf{O}^{\otimes}_{\neq i}\times
\prod(\mathbf{P}^{\otimes})^{\rm op}$
and $f_{(\mathbf{x},\mathbf{y}')}$
is an oplax $\mathcal{P}_j$-monoidal
functor for any
$(\mathbf{x},\mathbf{y}')\in
\prod\mathbf{O}^{\otimes}\times
\prod(\mathbf{P}^{\otimes}_{\neq j})^{\rm op}$.
\fi
 
%\newpage
\if
%Now, we consider another description
%of $\mathcal{O}_1$-$\mathcal{O}_2$-monoidal $\infty$-categoreis.
First, we recall mixed fibrations
introduced in \cite{Torii1}.
Let $p: X\to S\times T$ be a categorical fibration
where $S$ and $T$ are $\infty$-categories.
We say that $p: X\to S\times T$ is a mixed fibration
over $(S,T)$ if the following two conditions hold:
\begin{enumerate}
\item
The composite $\pi_S\circ p: X\to S$ is a coCartesian fibration
and the map $p$ carries $(\pi_S\circ p)$-coCartesian morphisms
to $\pi_S$-coCartesian morphisms,
where $\pi_S: S\times T\to S$ is the projection. 
\item
The composite $\pi_T\circ p: X\to T$ is a Cartesian fibration
and the map $p$ carries $(\pi_T\circ p)$-Cartesian morphisms
to $\pi_T$-Cartesian morphisms,
where $\pi_T: S\times T\to T$ is the projection. 
\end{enumerate}

Let $p: X\to S\times T$ be a mixed fibration
over $(S,T)$.
For $s\in S$ and $t\in T$,
we let 
$X_{s,\bullet}$ and $X_{\bullet,t}$
be the fibers of $\pi_S\circ p: X\to S$
and $\pi_T\circ p: X\to T$
at $s$ and $t$, respectively.
We denote by
$p_{s,\bullet}: X_{s,\bullet}\to T$ 
and $p_{\bullet,t}: X_{\bullet,t}\to S$
the restrictions of $p$.
Note that $p_{s,\bullet}: X_{s,\bullet}\to T$ 
is a Cartesian fibration for each $s\in S$,
and that $p_{\bullet,t}: X_{\bullet,t}\to S$
is a coCartesian fibration for each $t\in T$.

Next,
we will define $\mathcal{O}_1$-$\mathcal{O}_2$-monoidal
$\infty$-categories.
Let $p_1:\mathcal{O}_1^{\otimes}\to\Lambda(\Phi_1)$ and 
$p_2: \mathcal{O}_2^{\otimes}\to \Lambda(\Phi_2)$
be $\infty$-operads over 
perfect operator categories $\Phi_1$ and $\Phi_2$,
respectively.

\begin{definition}\rm
\label{definition:another-description-monoidal-by-mixed-fibration}
Let $p: \mathcal{C}^{\otimes}\to \mathcal{O}_1^{\otimes}\times
(\mathcal{O}_2^{\otimes})^{\rm op}$ be a mixed
fibration over 
$(\mathcal{O}_1^{\otimes},(\mathcal{O}_2^{\otimes})^{\rm op})$.
We say that $p: \mathcal{C}^{\otimes}\to \mathcal{O}_1^{\otimes}\times
(\mathcal{O}_2^{\otimes})^{\rm op}$ is 
an $\mathcal{O}_1$-$\mathcal{O}_2$-monoidal
$\infty$-category if
it satisfies the following two conditions:
\begin{enumerate}
\item
\label{condition1:o1-o2-monoidal}
The Segal morphism
\[ \mathcal{C}_{\bullet,y}^{\otimes}\longrightarrow 
   \prod_{i\in |I_2|}\!\!{}^{\mathcal{O}_1^{\otimes}}\
   \mathcal{C}_{\bullet,y_i}^{\otimes} \]
is an equivalence in $\cocart{\mathcal{O}_1^{\otimes}}$ 
%$\overcat{\mathcal{O}_1^{\otimes}}$
for each $y\in \mathcal{O}_2^{\otimes}$,
where 
$I_2=p_2(y)$,
$\{y\to y_i|\ i\in |I_2|\}$ are inert morphisms 
of $\mathcal{O}_2^{\otimes}$ lying over the inert morphisms
$\{I_2\to \{i\}|\ i\in |I_2|\}$,
and 
the right hand side is a product 
in $\cocart{\mathcal{O}_1^{\otimes}}$.

\item
\label{condition2:o1-o2-monoidal}
The Segal morphism
\[ \mathcal{C}_{x,\bullet}^{\otimes}\longrightarrow 
   \prod_{i\in |I_1|}\!\!{}^{(\mathcal{O}_2^{\otimes})^{\rm op}}\
   \mathcal{C}_{x_i,\bullet}^{\otimes}   \]
is an equivalence in $\cart{(\mathcal{O}_2^{\otimes})^{\rm op}}$ 
for each $x\in\mathcal{O}_1^{\otimes}$,
where $I_1=p_1(x)$,
$\{x\to x_i|\ i\in |I_1|\}$ are
inert morphisms of $\mathcal{O}_1^{\otimes}$
lying over the inert morphisms 
$\{I_1\to\{i\}|\ i\in |I_1|\}$,
and 
the right hand side is a product 
in $\cart{(\mathcal{O}_2^{\otimes})^{\rm op}}$.
\end{enumerate}
\end{definition}

\if0
Let $p: \mathcal{C}^{\otimes}\to \mathcal{O}_1^{\otimes}\times
(\mathcal{O}_2^{\otimes})^{\rm op}$ be an 
$\mathcal{O}_1$-$\mathcal{O}_2$-monoidal $\infty$-category.
For $X\in\mathcal{O}_1^{\otimes}$,
we denote by $p_{X,\bullet}: \mathcal{C}^{\otimes}_{X,\bullet}\to
(\mathcal{O}_2^{\otimes})^{\rm op}$
the composite $\pi_2\circ p$,
where $\pi_2: \mathcal{O}_1^{\otimes}\times 
(\mathcal{O}_2^{\otimes})^{\rm op}\to 
(\mathcal{O}_2^{\otimes})^{\rm op}$ is the projection.
Note that $p_{X,\bullet}: \mathcal{D}^{\otimes}_{X,\bullet}\to
(\mathcal{O}_2^{\otimes})^{\rm op}$
is a Cartesian fibration for any $X\in\mathcal{O}_1^{\otimes}$.
For $Y\in (\mathcal{O}_2^{\otimes})^{\rm op}$,
we also denote by $p_{\bullet,Y}: \mathcal{D}^{\otimes}_{\bullet,Y}\to
\mathcal{O}_1^{\otimes}$
the composite $\pi_1\circ p$,
where $\pi_1: \mathcal{O}_1^{\otimes}\times 
(\mathcal{O}_2^{\otimes})^{\rm op}\to 
\mathcal{O}_1^{\otimes}$ is the projection.
Note that $p_{\bullet,Y}: \mathcal{D}^{\otimes}_{\bullet,Y}\to
\mathcal{O}_1^{\otimes}$
is a coCartesian fibration for any 
$Y\in (\mathcal{O}_2^{\otimes})^{\rm op}$.
\fi

As in the case of duoidal $\infty$-categories,
we have three kinds of functors between
$\mathcal{O}_1$-$\mathcal{O}_2$-monoidal $\infty$-categories.
Now, we will define bilax monoidal functors between
$\mathcal{O}_1$-$\mathcal{O}_2$-monoidal $\infty$-categories.

\begin{definition}
\label{definition:bilax-monoida-functor}
\rm
Let $p: \mathcal{C}^{\otimes}\to\mathcal{O}_1^{\otimes}
\times (\mathcal{O}_2^{\otimes})^{\rm op}$
and $q: \mathcal{D}^{\otimes}\to\mathcal{O}_1^{\otimes}
\times (\mathcal{O}_2^{\otimes})^{\rm op}$
be $\mathcal{O}_1$-$\mathcal{O}_2$-monoidal
$\infty$-categories.
A bilax monoidal functor between $\mathcal{C}$ and $\mathcal{D}$
is a map $f: \mathcal{C}^{\otimes}\to \mathcal{D}^{\otimes}$
over $\mathcal{O}_1^{\otimes}\times (\mathcal{O}_2^{\otimes})^{\rm op}$
satisfying the following two conditions:

\begin{enumerate}

\item
\label{condition1:dual-inert-preservation}
$f$ carries $(\pi_1\circ p)$-coCartesian morphisms
over inert morphisms of $\mathcal{O}_1^{\otimes}$
to $(\pi_1\circ q)$-coCartesian morphisms,
where $\pi_1: \mathcal{O}_1^{\otimes}\times
(\mathcal{O}_2^{\otimes})^{\rm op}\to \mathcal{O}_1^{\otimes}$
is the projection.

\item
\label{condition2:dual-inert-preservation}
$f$ carries $(\pi_2\circ p)$-Cartesian morphisms
over inert morphisms of $(\mathcal{O}_2^{\otimes})^{\rm op}$
to $(\pi_2\circ q)$-Cartesian morphisms,
where $\pi_2: \mathcal{O}_1^{\otimes}\times
(\mathcal{O}_2^{\otimes})^{\rm op}\to (\mathcal{O}_2^{\otimes})^{\rm op}$
is the projection.

\end{enumerate}
\end{definition}

\begin{remark}\rm
We notice that the two conditions
in Definition~\ref{definition:bilax-monoida-functor}
are equivalent to
the following two conditions:
\begin{enumerate}

\item
For each $y\in (\mathcal{O}_2^{\otimes})^{\rm op}$,
the induced map
$f_{\bullet,y}: \mathcal{C}_{\bullet,y}^{\otimes}\to
\mathcal{D}_{\bullet,y}^{\otimes}$
carries $p_{\bullet,y}$-coCartesian morphisms over 
inert morphisms of $\mathcal{O}_1^{\otimes}$
to $q_{\bullet,y}$-coCartesian morphisms.

\item
For each $x\in \mathcal{O}_1^{\otimes}$,
the induced map
$f_{x,\bullet}: \mathcal{C}_{x,\bullet}^{\otimes}\to
\mathcal{D}_{x,\bullet}^{\otimes}$
carries $p_{x,\bullet}$-Cartesian morphisms over 
inert morphisms of $(\mathcal{O}_2^{\otimes})^{\rm op}$
to $q_{x,\bullet}$-Cartesian morphisms.

\end{enumerate}

\end{remark}

\if0
For an $\infty$-operad $p: \mathcal{O}^{\otimes}\to \Lambda(\Phi)$
over a perfect operator category $\Phi$, 
we let $\mathcal{O}^{\otimes,\natural}$ be the 
simplicial set $\mathcal{O}^{\otimes}$
equipped with the set of $p$-coCartesian morphisms
over inert morphisms of $\Lambda(\Phi)$
as marked edges.
%Let $\mathcal{O}_1^{\otimes}$ and
%$\mathcal{O}_2^{\otimes}$ be
%$\infty$-operads over $\Phi_1$ and $\Phi_2$, respectively.
\fi

\begin{definition}\rm
We write 
\[ \mathcal{O}_1\mbox{\rm -}\mathcal{O}_2\mbox{\rm -}
   {\rm Mon}^{\rm bilax} \]
for the subcategory
of $\overcat{\mathcal{O}_1^{\otimes}\times
(\mathcal{O}_2^{\otimes})^{\rm op}}$
with $\mathcal{O}_1$-$\mathcal{O}_2$-monoidal
$\infty$-categories as objects and
bilax monoidal functors as morphisms.
\end{definition}

%\bigskip

As in \cite[\S\S4.2 and 4.3]{Torii1},
we can describe 

%\begin{proposition}
\begin{theorem}
\if0
Let $\mathcal{O}_1^{\otimes}$ and
$\mathcal{O}_2^{\otimes}$ be
$\infty$-operads over perfect operator categories
$\Phi_1$ and $\Phi_2$, respectively.
\fi
There are equivalences
\[ \begin{array}{rcl}
    \mathcal{O}_1\mbox{-}\mathcal{O}_2\mbox{-}\duo^{\rm bilax}
    &\simeq&
    \laxmon_{\mathcal{O}_1}(\oplaxmoncatp{\mathcal{O}_2})\\[2mm]
    &\simeq&
    \oplaxmon_{\mathcal{O}_2}(\laxmoncatp{\mathcal{O}_1})\\
   \end{array}\]
of $\infty$-categories.
\end{theorem}
%\end{proposition}

Next, as in \cite[\S\S4.2 and 4.3]{Torii1},
we will show that
$\mathcal{O}_1$-$\mathcal{O}_2$-monoidal $\infty$-categories
and bilax monoidal functors
can be described in terms of 
$\laxmoncatp{\mathcal{O}_1}$
and $\oplaxmoncatp{\mathcal{O}_2}$.
Since the $\infty$-category
$\laxmoncatp{\mathcal{O}_1}$  has finite products,
we can consider
$\mathcal{O}_1$-monoid objects
of $\oplaxmoncatp{\mathcal{O}_2}$.
Also, 
since the $\infty$-category
$\oplaxmoncatp{\mathcal{O}_2}$ has finite products,
we can consider
$\mathcal{O}_2$-monoid objects
of $\laxmoncatp{\mathcal{O}_1}$.

\begin{proposition}
\label{prop:idetifying-O1-O2-monoidal-to-monoid-oplaxmonodal}
We can identify
an $\mathcal{O}_1$-$\mathcal{O}_2$-monoidal $\infty$-category
with 
an $\mathcal{O}_1$-monoid object of
$\oplaxmoncatp{\mathcal{O}_2}$.
Also,
we can identify
an $\mathcal{O}_1$-$\mathcal{O}_2$-monoidal $\infty$-category
with 
an $\mathcal{O}_2$-monoid object of
$\laxmoncatp{\mathcal{O}_1}$.
\end{proposition}

\proof
We shall prove that 
an $\mathcal{O}_1$-$\mathcal{O}_2$-monoidal $\infty$-category
can be identified with 
an $\mathcal{O}_1$-monoid object of
$\oplaxmoncatp{\mathcal{O}_2}$.
The other case can be shown in the similar way.

Suppose that $p: \mathcal{C}^{\otimes}\to
\mathcal{O}_1^{\otimes}\times(\mathcal{O}_2^{\otimes})^{\rm op}$
is a $\mathcal{O}_1$-$\mathcal{O}_2$-monoidal 
$\infty$-category.
Since $p$ is a mixed fibration,
we have an associated functor
$C: \mathcal{O}_1^{\otimes}\to 
(\cart{(\mathcal{O}_2^{\otimes})^{\rm op}})_{\rm oplax}$.
Conditions (\ref{condition2:o1-o2-monoidal}) 
and
(\ref{condition2:o1-o2-monoidal}) in 
Definition~\ref{definition:another-description-monoidal-by-mixed-fibration}
implies that
$C$ factors through $\oplaxmoncatp{\mathcal{O}_2}$.
Furthermore,
Condition (\ref{condition2:o1-o2-monoidal}) in 
Definition~\ref{definition:another-description-monoidal-by-mixed-fibration}
ensures that $C$
is an $\mathcal{O}_1$-monoid object
of $\oplaxmoncatp{\mathcal{O}_2}$.  

Conversely,
suppose that $\mathcal{D}: \mathcal{O}_1^{\otimes}
\to \oplaxmoncatp{\mathcal{O}_2}$
is an $\mathcal{O}_1$-monoid object of
$\oplaxmoncatp{\mathcal{O}_2}$.
The composite 
$\mathcal{O}_1^{\otimes}\stackrel{\mathcal{D}}{\to}
\oplaxmoncatp{\mathcal{O}_2}\to
\cart{(\mathcal{O}_2^{\otimes})^{\rm op}}$
corresponds to
a mixed fibration
$p: \mathcal{D}^{\otimes}\to \mathcal{O}_1^{\otimes}
\times (\mathcal{O}_2^{\otimes})^{\rm op}$
by {\color{red} \cite[Proposition~3.28]{Torii1}}.
For any morphism $X\to X'$ in $\mathcal{O}_1^{\otimes}$
and any object $Y\in (\mathcal{O}_2^{\otimes})^{\rm op}$,
since $\mathcal{D}^{\otimes}_{X,\bullet}\to
\mathcal{D}^{\otimes}_{X',\bullet}$
is a morphism of $\oplaxmoncatp{\mathcal{O}_2}$,
we have a commutative diagram
\[ \begin{array}{ccc}
    \mathcal{D}^{\otimes}_{X,Y}&\longrightarrow&
    \prod_{i\in |I_2|}\mathcal{D}^{\otimes}_{X,Y_i}\\[2mm]
    \bigg\downarrow & & \bigg\downarrow \\[4mm]
    \mathcal{D}^{\otimes}_{X',Y}&\longrightarrow& 
    \prod_{i\in |I_2|}\mathcal{D}^{\otimes}_{X',Y_i},\\
   \end{array} \]
where $I_2=p_2(Y)$ and 
the horizontal arrows are induced by
Cartesian morphisms $\{Y\to Y_i|\ i\in |I_2|\}$
lying over inert morphisms $\{I_2\to\{i\}|\ i\in |I_2|\}$.
This implies that 
the Segal map 
$\mathcal{D}^{\otimes}_{\bullet,Y}\to \prod_{i\in |I_2|}^{\mathcal{O}_1^{\otimes}}
\mathcal{D}^{\otimes}_{\bullet,Y_i}$
preserves coCartesian morphisms.
Since the horizontal arrows of the above commutative diagram
are equivalences,
the Segal map 
$\mathcal{D}^{\otimes}_{\bullet,Y}\to \prod_{i\in |I_2|}^{\mathcal{O}_1^{\otimes}}
\mathcal{D}^{\otimes}_{\bullet,Y_i}$
is an equivalence in $\cocart{\mathcal{O}_1^{\otimes}}$.
Furthermore,
since $\mathcal{D}$ is an $\mathcal{O}_1$-monoid object
of $\oplaxmoncatp{\mathcal{O}_2}$,
the Segal map
$\mathcal{D}^{\otimes}_{X,\bullet}\to 
\prod_{i\in |I_1|}^{(\mathcal{O}_2^{\otimes})^{\rm op}}
\mathcal{D}^{\otimes}_{X_i,\bullet}$
is an equivalence in $\cart{(\mathcal{O}_2^{\otimes})^{\rm op}}$
for any $X\in\mathcal{O}_1^{\otimes}$.
This completes the proof.
\qed

%\begin{proposition}
\begin{theorem}
\if0
Let $\mathcal{O}_1^{\otimes}$ and
$\mathcal{O}_2^{\otimes}$ be
$\infty$-operads over perfect operator categories
$\Phi_1$ and $\Phi_2$, respectively.
\fi
There are equivalences
\[ \begin{array}{rcl}
    \mathcal{O}_1\mbox{-}\mathcal{O}_2\mbox{-}\duo^{\rm bilax}
    &\simeq&
    \laxmon_{\mathcal{O}_1}(\oplaxmoncatp{\mathcal{O}_2})\\[2mm]
    &\simeq&
    \oplaxmon_{\mathcal{O}_2}(\laxmoncatp{\mathcal{O}_1})\\
   \end{array}\]
of $\infty$-categories.
\end{theorem}
%\end{proposition}

\proof
We shall prove an equivalence
$\mathcal{O}_1\mbox{-}\mathcal{O}_2\mbox{-}\duo^{\rm bilax}
 \simeq
 \laxmon_{\mathcal{O}_1}(\oplaxmoncatp{\mathcal{O}_2})$.
By Proposition~\ref{prop:idetifying-O1-O2-monoidal-to-monoid-oplaxmonodal},
we can identify the objects of
$\laxmon_{\mathcal{O}_1}(\oplaxmoncatp{\mathcal{O}_2})$
with $\mathcal{O}_1$-$\mathcal{O}_2$-monoidal $\infty$-categories.
We can verify that
morphisms of $\laxmon_{\mathcal{O}_1}(\oplaxmoncatp{\mathcal{O}_2})$
coincide with
bilax functors between
$\mathcal{O}_1$-$\mathcal{O}_2$-monoidal $\infty$-categories.
The other case can be shown in the similar way.
\qed

\begin{remark}\rm
Recall that $R: \cat\to\cat$
is the functor which assigns to an
$\infty$-category its opposite $\infty$-category.
If $p: \mathcal{D}^{\otimes}\to 
\mathcal{O}_1^{\otimes}\times(\mathcal{O}_2^{\otimes})^{\rm op}$
is an $\mathcal{O}_1$-$\mathcal{O}_2$-monoidal 
$\infty$-category,
then
$R(p): (\mathcal{D}^{\otimes})^{\rm op}\to
\mathcal{O}_2^{\otimes}\times (\mathcal{O}_1^{\otimes})^{\rm op}$
is an $\mathcal{O}_2$-$\mathcal{O}_1$-monoidal
$\infty$-category.
\end{remark}
\fi

%\newpage
%\input{adjoint}
\section{Adjoints of (op)lax monoidal functors}
\label{section:adjoints-monoidal-functors}

In this section we study 
adjoints of (op)lax monoidal functors
between higher monoidal $\infty$-categories.
For this purpose,
we study adjoints of (op)lax morphisms between
$\mathbf{S}$-(co)Cartesian fibrations
in \S\ref{subsection:adjoints-mixed-fibrations}.
In \S\ref{subsection:monoidal-adjoints}
we prove the main theorem
(Theorem~\ref{them:adjoint-monoidal-functors})
which says that 
the $\infty$-category
of coCartesian $\mathbf{O}$-monoidal
$\infty$-categories and right adjoint
lax $\mathbf{O}$-monoidal functors
is equivalent to the opposite of 
the $\infty$-category of Cartesian 
$\mathbf{O}_{\rm rev}$-monoidal
$\infty$-categories and left adjoint oplax 
$\mathbf{O}_{\rm rev}$-monoidal functors.

\subsection{Adjoints of (op)lax morphisms
between $\mathbf{S}$-(co)Cartesian fibrations}
\label{subsection:adjoints-mixed-fibrations}

In this subsection
we study adjoints of (op)lax morphisms 
between $\mathbf{S}$-(co)Cartesian fibrations.
We show that the $\infty$-category
of coCartesian $\mathbf{S}$-fibrations
and right adjoint lax $\mathbf{S}$-morphisms
is equivalent to the opposite of
the $\infty$-category of
Cartesian $\mathbf{S}^{\rm op}_{\rm rev}$-fibrations and
left adjoint oplax $\mathbf{S}^{\rm op}_{\rm rev}$-morphisms.

\begin{definition}\rm
Let $\mathbf{S}=(S_1,\ldots,S_n)$
be a finite sequence of $\infty$-categories.
We write
\[ \rlaxcocart{\mathbf{S}}(\cat) \]
for the wide subcategory of 
$\laxcocart{\mathbf{S}}(\cat)$
spanned by those lax $\mathbf{S}$-morphisms $f: X\to Y$
over $\prod \mathbf{S}$
that satisfy the following condition:
%\begin{enumerate}
%\item[]
For each $\mathbf{s}%=(s_1,\ldots,s_n)
\in \prod\mathbf{S}$,
the restriction
$f_{\mathbf{s}}: X_{\mathbf{s}}\to Y_{\mathbf{s}}$
is a right adjoint functor.
%\end{enumerate}

Dually, we write
\[ \loplaxcart{\mathbf{S}}(\cat) \]
for the wide subcategory of 
$\oplaxcart{\mathbf{S}}(\cat)$
spanned by those oplax $\mathbf{S}$-morphisms $f: X\to Y$
over $\prod \mathbf{S}$
that satisfy the following condition:
%\begin{enumerate}
%\item[]
For each $\mathbf{s}%=(s_1,\ldots,s_n)
\in \prod\mathbf{S}$,
the restriction
$f_{\mathbf{s}}: X_{\mathbf{s}}\to Y_{\mathbf{s}}$
is a left adjoint functor.
%\end{enumerate}
\end{definition}

The goal of this subsection 
is to prove the following theorem. 

\begin{theorem}
\label{thm:equivalanece-higher-adjoints}
There is a natural equivalence
\[ \rlaxcocart{\mathbf{S}}(\cat)\simeq
   \loplaxcart{\mathbf{S}^{\rm op}_{\rm rev}}(\cat)^{\rm op}\]
of $\infty$-categories,
which is given 
on objects by taking dual fibrations 
and on morphisms by taking adjoints fiberwise.
\end{theorem}

We recall that a biCartesian fibration
is both a Cartesian fibration and a coCartesian fibration.
See \cite[\S4.7.4]{Lurie2} for the relationship between
biCartesian fibrations and adjoint functors.
In particular, if
a coCartesian fibration $p: X\to S$
in which the induced functor $X_s\to X_{s'}$
admits a right adjoint for each morphism
$s\to s'$ in S,
then $p$ is a biCartesian morphism.
Dually, if a Cartesian fibration $q: Y\to T$
in which the induced functor $Y_t\to Y_{t'}$
admits a left adjoint for each morphism
$t'\to t$ in $T$,
then $q$ is a biCartesian fibration.

\begin{notation}\rm
We write 
\[ {\rm bCart}_{/U}^{\rm bilax}(\cat) \]
for the full subcategory of
$\overcat{U}$ spanned by
biCartesian fibrations over $U$.
\end{notation}

\if0
\begin{proposition}
There are natural equivalences
\[ {\rm Mix}_{/(\mathbf{S},U)}^{\rm oplax}
   (\cat)^{\mbox{\scriptsize $\mathbf{S}$-R}}\simeq
   {\rm Fun}\left(\prod\mathbf{S}, {\rm bCart}_{/U}^{\rm bilax}
   (\cat)\right)
   \simeq
   {\rm Mix}_{/(U,\mathbf{S}^{\rm op})}^{\rm lax}
   (\cat)^{\mbox{\scriptsize $\mathbf{S}^{\rm op}$-L}}\]
of $\infty$-categories.
\end{proposition}
\fi

%We can easily verify the following proposition.

\begin{lemma}
\label{lemma:mixed-adjoint-bicartesian}
There are natural equivalences
\[ \begin{array}{rcl}
%   {\rm Mix}_{/(\mathbf{S},U)}^{\rm bilax}
%   (\cat)^{\mbox{\scriptsize $\mathbf{S}$-R}}
   \oplaxcart{U^{\rm op}}(\rlaxcocart{\mathbf{S}}(\cat))
   &\simeq&
   \laxcocart{\mathbf{S}}
   ({\rm bCart}_{/U}^{\rm bilax}(\cat)),\\[2mm]
   \laxcocart{U}(\loplaxcart{\mathbf{S}}(\cat))
%   {\rm Mix}_{/(U,\mathbf{T})}^{\rm bilax}
%   (\cat)^{\mbox{\scriptsize $\mathbf{T}$-L}}
   &\simeq&
   \oplaxcart{\mathbf{S}}
   ({\rm bCart}_{/U}^{\rm bilax}(\cat))\\
   \end{array}\]
of $\infty$-categories.
\end{lemma}

\proof
The first equivalence follows from the fact
that the both sides are full subcategories
of ${\rm Mix}_{/(\mathbf{S},U)}^{\rm bilax}(\cat)$
spanned by those mixed fibrations
$p: X\to \prod\mathbf{S}\times U$
such that $p_{\mathbf{s}}: X_{\mathbf{s}}\to U$
is a coCartesian fibration for each $\mathbf{s}\in\prod\mathbf{S}$.
The second equivalence can be proved similarly.
\qed

\proof[Proof of Theorem~\ref{thm:equivalanece-higher-adjoints}]
%By the Yoneda Lemma,
%it suffices to show that
First, we shall show that  
there exists a natural equivalence
\[ {\rm Map}_{\cat}(U,\rlaxcocart{\mathbf{S}}(\cat))
    \simeq
    {\rm Map}_{\cat}(U,\loplaxcart{\mathbf{S}^{\rm op}_{\rm rev}}   
    (\cat)^{\rm op}) \]
of $\infty$-groupoids 
for any $\infty$-category $U$.

By Lemmas~\ref{lemma:straightening-unstraightening-Z} 
and \ref{lemma:strict-oplax-tilde-equivalence},
we have natural equivalences
\[ \begin{array}{rcl}
    {\rm Map}_{\cat}(U,\rlaxcocart{\mathbf{S}}(\cat))
    &\simeq&
    \cart{U^{\rm op}}(\rlaxcocart{\mathbf{S}}(\cat))^{\simeq}\\[2mm]
    &\simeq&
    \oplaxcart{U^{\rm op}}(\rlaxcocart{\mathbf{S}}(\cat))^{\simeq}\\%[2mm]
%    &\simeq&
%    {\rm Mix}_{/(\mathbf{S},U^{\rm op})}^{\rm bilax}
%     (\cat)^{\mbox{\scriptsize $\mathbf{S}$-R},\simeq}.\\
   \end{array}\]
Similarly,
by Lemmas~\ref{lemma:Z-(un)straightening} 
and \ref{lemma:strict-lax-tilde-equivalence},
we have natural equivalences
\[ \begin{array}{rcl}
    {\rm Map}_{\cat}(U,\loplaxcart{\mathbf{S}^{\rm op}_{\rm rev}}   
    (\cat)^{\rm op})
    &\simeq&
    {\rm Map}_{\cat}(U^{\rm op},\loplaxcart{\mathbf{S}^{\rm op}_{\rm rev}}
    (\cat))\\[2mm]
    &\simeq&
    \cocart{U^{\rm op}}
    (\loplaxcart{\mathbf{S}^{\rm op}_{\rm rev}}(\cat))^{\simeq}\\[2mm]
    &\simeq&
    \laxcocart{U^{\rm op}}
    (\loplaxcart{\mathbf{S}^{\rm op}_{\rm rev}}(\cat))^{\simeq}\\%[2mm]
%    &\simeq&
%    {\rm Mix}_{/(U^{\rm op},\mathbf{S}^{\rm op}_{\rm rev})}^{\rm bilax}
%     (\cat)^{\mbox{\scriptsize $\mathbf{S}^{\rm op}_{\rm rev}$-L},\simeq}.\\
   \end{array}\]
By Lemma~\ref{lemma:mixed-adjoint-bicartesian}
and Corollary~\ref{cor:cocart-cart-duality},
\if0
there are natural equivalences
\[ \begin{array}{rcl}
    {\rm Mix}_{/(\mathbf{S},U^{\rm op})}^{\rm bilax}
    (\cat)^{\mbox{\scriptsize $\mathbf{S}$-R},\simeq}
     &\simeq&
   \laxcocart{\mathbf{S}}
   ({\rm bCart}_{/U^{\rm op}}^{\rm bilax}(\cat))^{\simeq}\\[2mm]
    &\simeq&
   \oplaxcart{\mathbf{S}^{\rm op}_{\rm rev}}
    ({\rm bCart}_{/U^{\rm op}}^{\rm bilax}
    (\cat))^{\simeq}\\[2mm]
    &\simeq&
    {\rm Mix}_{/(U^{\rm op},\mathbf{S}^{\rm op}_{\rm rev})}^{\rm bilax}
     (\cat)^{\mbox{\scriptsize $\mathbf{S}^{\rm op}_{\rm rev}$-L},\simeq}\\[2mm]
   \end{array}\]
\fi
%Thus, 
we obtain the desired equivalence.

By the Yoneda Lemma,
we obtain an equivalence between
the $\infty$-categories
$\rlaxcocart{\mathbf{S}}(\cat)$
and
$\loplaxcart{\mathbf{S}^{\rm op}_{\rm rev}}(\cat)^{\rm op}$.
We can verify that this equivalence
is given on objects by taking dual fibrations
and on morphisms by taking adjoints fiberwise.
\if0
We have natural equivalences
\[ \begin{array}{rcl}
    {\rm Map}_{\cat}(U,\rlaxcocart{\mathbf{S}}(\cat))
    &\simeq&
    \cart{U^{\rm op}}(\rlaxcocart{\mathbf{S}}(\cat))^{\simeq}\\[2mm]
    &\simeq&
    \oplaxcart{U^{\rm op}}(\rlaxcocart{\mathbf{S}}(\cat))^{\simeq}\\[2mm]
    &\simeq&
    {\rm Mix}_{/(\mathbf{S},U^{\rm op})}^{\rm bilax}
     (\cat)^{\mbox{\scriptsize $\mathbf{S}$-R},\simeq}\\[2mm]
    &\simeq&
    \laxcocart{\mathbf{S}}({\rm bCart}_{/U^{\rm op}}^{\rm bilax}
    (\cat))^{\simeq}\\[2mm]
    &\simeq&
   \oplaxcart{\mathbf{S}^{\rm op}_{\rm rev}}
    ({\rm bCart}_{/U^{\rm op}}^{\rm bilax}
    (\cat))^{\simeq}\\[2mm]
    &\simeq&
    {\rm Mix}_{/(U^{\rm op},\mathbf{S}^{\rm op}_{\rm rev})}^{\rm bilax}
     (\cat)^{\mbox{\scriptsize $\mathbf{S}^{\rm op}_{\rm rev}$-L},\simeq}\\[2mm]
    &\simeq&
    \laxcocart{U^{\rm op}}
    (\loplaxcart{\mathbf{S}^{\rm op}_{\rm rev}}(\cat))^{\simeq}\\[2mm]
    &\simeq&
    \cocart{U^{\rm op}}
    (\loplaxcart{\mathbf{S}^{\rm op}_{\rm rev}}(\cat))^{\simeq}\\[2mm]
    &\simeq&
    {\rm Map}_{\cat}(U^{\rm op},\loplaxcart{\mathbf{S}^{\rm op}_{\rm rev}}
    (\cat))\\[2mm]
    &\simeq&
    {\rm Map}_{\cat}(U,\loplaxcart{\mathbf{S}^{\rm op}_{\rm rev}}   
    (\cat)^{\rm op})\\
   \end{array}\]
of $\infty$-groupoids.
\fi
\qed

%\newpage
%\input{monoidal-adjoint}
\subsection{Monoidal adjoints of higher monoidal $\infty$-categories}
\label{subsection:monoidal-adjoints}

Let $\mathbf{O}^{\otimes}=(\mathcal{O}_1^{\otimes},
\ldots,\mathcal{O}_n^{\otimes})$
be a finite sequence of $\infty$-operads
over perfect operator categories.
In this subsection
we show that the $\infty$-category
of coCartesian $\mathbf{O}$-monoidal
$\infty$-categories and right adjoint 
lax $\mathbf{O}$-monoidal functors
is equivalent to
the opposite of the $\infty$-category of  
Cartesian $\mathbf{O}_{\rm rev}$-monoidal 
$\infty$-categories
and left adjoint oplax $\mathbf{O}_{\rm rev}$-monoidal
functors
(Theorem~\ref{them:adjoint-monoidal-functors}).

\begin{definition}\rm
We define an $\infty$-category
\[ \mathsf{Mon}_{\mathbf{O}}^{\rm lax,R}(\cat) \]
to be the wide subcategory of 
$\mathsf{Mon}_{\mathbf{O}}^{\rm lax}(\cat)$
spanned by those lax $\mathbf{O}$-monoidal
functors $f: \mathcal{C}^{\otimes}\to \mathcal{D}^{\otimes}$
that satisfy the following condition:
%
%\begin{enumerate}
%\item[]
%For each $x=(x_1,\ldots,x_n)\in \mathcal{O}_1\times\cdots\times\mathcal{O}_n$,
For each $x\in\prod\mathbf{O}$,
the restriction
$f_x: \mathcal{C}_x\to \mathcal{D}_x$
is a right adjoint functor.
%\end{enumerate}

Dually,
we define an $\infty$-category
\[ \mathsf{Mon}_{\mathbf{O}}^{\rm oplax,L}(\cat) \]
to be the wide subcategory of 
$\mathsf{Mon}_{\mathbf{O}}^{\rm oplax}(\cat)$
spanned by those oplax $\mathbf{O}$-monoidal
functors $f: \mathcal{C}^{\otimes}\to 
\mathcal{D}^{\otimes}$
that satisfy the following condition:
%
%\begin{enumerate}
%\item[]
%For each $x=(x_1,\ldots,x_n)\in
%\mathcal{O}_1\times\cdots\times\mathcal{O}_n$,
For each $x\in\prod\mathbf{O}^{\rm op}$,
the restriction
$f_x: \mathcal{C}_x\to \mathcal{D}_x$
is a left adjoint functor.
%\end{enumerate}
\end{definition}

\begin{remark}\rm
Notice that 
$\mathsf{Mon}_{\mathbf{O}}^{\rm lax,R}(\cat)$
is a subcategory of 
$\rlaxcocart{\mathbf{O}^{\otimes}}(\cat)$.
This follows from the fact that
the restriction $f_{\mathbf{x}}: \mathcal{C}^{\otimes}_{\mathbf{x}}
\to \mathcal{D}^{\otimes}_{\mathbf{x}}$
for $\mathbf{x}=(x_1,\ldots,x_n)\in \prod\mathbf{O}^{\otimes}$
with $p(x_i)=I_i\ (1\le i\le n)$
is equivalent to a product
of $f_x$ for $x=(x_{1,i_1},\ldots, x_{n,i_n})
\in \prod\mathbf{O}$
over $i_1\in |I_1|,\ldots,i_n\in |I_n|$.
Similarly,
$\mathsf{Mon}_{\mathbf{O}}^{\rm oplax,L}(\cat)$
is a subcategory of 
$\loplaxcart{(\mathbf{O}^{\otimes})^{\rm op}}(\cat)$.
\end{remark}

The following is 
the main theorem of this paper.
 
\begin{theorem}
\label{them:adjoint-monoidal-functors}
The equivalence of
Theorem~\ref{thm:equivalanece-higher-adjoints}
restricts to an equivalence
\[ \mathsf{Mon}_{\mathbf{O}}^{\rm lax,R}(\cat)
   \simeq
   \mathsf{Mon}_{\mathbf{O}_{\rm rev}}^{\rm oplax,L}(\cat)^{\rm op} \]
of $\infty$-categories,
which is given on objects by taking
dual fibrations and on morphisms by taking adjoints fiberwise. 
\end{theorem}

\begin{corollary}
The left adjoint of a lax $\mathbf{O}$-monoidal
functor between coCartesian $\mathbf{O}$-monoidal 
$\infty$-categories
is canonically an oplax $\mathbf{O}_{\rm rev}$-monoidal 
functor between the corresponding
Cartesian $\mathbf{O}_{\rm rev}$-monoidal $\infty$-categories, 
and vice versa.
\end{corollary}

\if0
In order to prove Theorem~\ref{them:adjoint-monoidal-functors},
we need the following lemma.
{\color{red} 2021/10/07 前のsectionで示しているか確認。
2021/10/12
Lemma~{lemma:lax-oplax-equivalence-groupoids}
と同じ}

\begin{lemma}[Lemma~\ref{lemma:lax-oplax-equivalence-groupoids}]
\label{lemma:laxmon-oplaxmon-groupoid-commutes}
Let $\mathcal{Z}$ be a subcategory of $\overcat{Z}$
that is closed under finite products and equivalences.
The equivalence
$\laxcocart{\mathcal{O}^{\otimes}}(\mathcal{Z})^{\simeq}\simeq
\oplaxcart{\mathcal{O}^{\otimes,{\rm op}}}
(\mathcal{Z})^{\simeq}$
restricts to an equivalence
\[ \mathsf{Mon}_{\mathcal{O}}^{\rm lax}(\mathcal{Z})^{\simeq}
   \simeq
   \mathsf{Mon}_{\mathcal{O}}^{\rm oplax}(\mathcal{Z})^{\simeq} \]
of $\infty$-groupoids.
\end{lemma}

\proof
The equivalence
$\laxcocart{\mathcal{O}^{\otimes}}(\mathcal{Z})^{\simeq}\simeq
\oplaxcart{\mathcal{O}^{\otimes,{\rm op}}}
(\mathcal{Z})^{\simeq}$
is given by the composite of the following equivalences
\[ \laxcocart{\mathcal{O}^{\otimes}}(\mathcal{Z})^{\simeq}\simeq
   {\rm Map}_{\cat}(\mathcal{O}^{\otimes},\mathcal{Z})
   \simeq
   \oplaxcart{\mathcal{O}^{\otimes,{\rm op}}}
   (\mathcal{Z})^{\simeq}.\]
The lemma follows from the following
equivalences
\[ \mathsf{Mon}_{\mathcal{O}}^{\rm lax}(\mathcal{Z})^{\simeq}
   \simeq
   {\rm Mon}_{\mathcal{O}}(\mathcal{Z})^{\simeq}
   \simeq
   \mathsf{Mon}_{\mathcal{O}}^{\rm oplax}(\mathcal{Z})^{\simeq}, \]
where ${\rm Mon}_{\mathcal{O}}(\mathcal{Z})$
is the $\infty$-category of monoid objects in $\mathcal{Z}$.
\qed
\fi

\proof[Proof of Theorem~\ref{them:adjoint-monoidal-functors}]
We will prove the theorem by induction 
on $l(\mathbf{O})$.
When $l(\mathbf{O})=0$,
this follows from Theorem~\ref{thm:equivalanece-higher-adjoints}
for $\mathbf{S}=\emptyset$.

Suppose $l(\mathbf{O})>0$.
First, we note that
an equivalence
\[ \rlaxcocart{\mathbf{S}}(\cat)^{\simeq}\simeq
   \loplaxcart{\mathbf{S}^{\rm op}_{\rm rev}}(\cat)^{\simeq}\]
induced by the equivalence
in Theorem~\ref{thm:equivalanece-higher-adjoints}
restricts to an equivalence
\[ \mathsf{Mon}_{\mathbf{O}}^{\rm lax,R}(\cat)^{\simeq}
   \simeq
   \mathsf{Mon}_{\mathbf{O}_{\rm rev}}^{\rm oplax,L}(\cat)^{\simeq} \]
of $\infty$-groupoids.

\if0
we will prove an equivalence 
\[ \mathsf{Mon}_{\mathbf{O}}^{\rm lax,R}(\cat)^{\simeq}
   \simeq
   \mathsf{Mon}_{\mathbf{O}_{\rm rev}}^{\rm oplax,L}(\cat)^{\simeq} \]
of $\infty$-groupoids.
Notice that There are equivalences
\[ \begin{array}{rcl}
     \mathsf{Mon}_{\mathbf{O}}^{\rm lax,R}(\cat)^{\simeq}
     &\simeq&
     \mathsf{Mon}_{\mathbf{O}}^{\rm lax}(\cat)^{\simeq},\\[2mm]
     \mathsf{Mon}_{\mathbf{O}_{\rm rev}}^{\rm oplax,L}(\cat)^{\simeq}
     &\simeq&
     \mathsf{Mon}_{\mathbf{O}_{\rm rev}}^{\rm oplax}(\cat)^{\simeq}
   \end{array} \]
of $\infty$-groupoids.
By using Lemma~\ref{lemma:laxmon-oplaxmon-groupoid-commutes},
we have equivalences
\[ \begin{array}{rcl}
    \mathsf{Mon}_{\mathbf{O}}^{\rm lax}(\cat)^{\simeq}
    &\simeq&
    \mathsf{Mon}_{\mathcal{O}_1}^{\rm lax}
    (\mathsf{Mon}_{\mathbf{O}_{\ge 2}}^{\rm lax}(\cat))^{\simeq}\\[2mm]
    &\simeq&
    \mathsf{Mon}_{\mathcal{O}_1}^{\rm oplax}
    (\mathsf{Mon}_{\mathbf{O}_{\ge 2}}^{\rm lax}(\cat))^{\simeq}\\[2mm]
    &\simeq&
    \mathsf{Mon}_{\mathbf{O}_{\ge 2}}^{\rm lax}
    (\mathsf{Mon}_{\mathcal{O}_1}^{\rm oplax}(\cat))^{\simeq}\\[2mm]
    &\simeq&
    \mathsf{Mon}_{(\mathbf{O}_{\ge 2})_{\rm rev}}^{\rm oplax}
    (\mathsf{Mon}_{\mathcal{O}_1}^{\rm oplax}(\cat))^{\simeq}\\[2mm]
    &\simeq&
    \mathsf{Mon}_{\mathbf{O}_{\rm rev}}^{\rm oplax}(\cat)^{\simeq}\\
   \end{array}\]
of $\infty$-groupoids.
Thus, we obtain an equivalence
$\mathsf{Mon}_{\mathbf{O}}^{\rm lax,R}(\cat)^{\simeq}
   \simeq
   \mathsf{Mon}_{\mathbf{O}_{\rm rev}}^{\rm oplax,L}(\cat)^{\simeq}$
of $\infty$-groupoids.
\fi

Thus, by symmetry,
it suffices to show the following claim:
Suppose that $f: \mathcal{C}^{\otimes}\to \mathcal{D}^{\otimes}$
is a lax $\mathbf{O}$-monoidal functor
between coCartesian $\mathbf{O}$-monoidal
$\infty$-categories.
We assume that $f_{\mathbf{x}}$ is a right adjoint
functor for each $\mathbf{x}\in \prod\mathbf{O}^{\otimes}$.
If $g: (\mathcal{D}^{\otimes})^{\vee}\to
(\mathcal{C}^{\otimes})^{\vee}$  
is a corresponding oplax
$(\mathbf{O}^{\otimes}_{\rm rev})^{\rm op}$-morphism
under the equivalence in Theorem~\ref{thm:equivalanece-higher-adjoints},
then $g$ is an oplax $\mathbf{O}_{\rm rev}$-monoidal functor.

\if0
Let $f: \mathcal{C}^{\otimes}\to\mathcal{D}^{\otimes}$
be a lax $\mathbf{O}$-monoidal functor.
We assume that $f_{\mathbf{x}}$ is a right adjoint
functor for each $\mathbf{x}\in \prod\mathbf{O}^{\otimes}$.
By Theorem~\ref{thm:equivalanece-higher-adjoints},
there is canonically
an oplax $\mathbf{O}^{\rm op}_{\rm rev}$-morphism
$g: \mathcal{D}^{\otimes,\vee}\to \mathcal{C}^{\otimes,\vee}$
such that $g_\mathbf{x}$ is a left adjoint
functor to $f_{\mathbf{x}}$ for each 
$\mathbf{x}\in\prod\mathbf{O}^{\otimes}$.
We have to show that $g$ is an oplax $\mathbf{O}$-monoidal
functor.
\fi

By Proposition~\ref{prop:criterion-(op)lax-monoidal-functor},
it suffices to show that
$g_{\mathbf{x}}: 
(\mathcal{D}^{\otimes})^{\vee}_{\mathbf{x}}\to
(\mathcal{C}^{\otimes})^{\vee}_{\mathbf{x}}$
is an oplax $\mathcal{O}_i$-monoidal
functor for each
$\mathbf{x}\in \prod\mathbf{O}^{\otimes}_{\neq i}
\ (1\le i\le l(\mathbf{O}))$.
Let $x\to x'$ be an inert morphism 
of $\mathcal{O}_i^{\otimes}$.
Since $f$ is a lax $\mathbf{O}$-monoidal functor,
$f_{\mathbf{x}}:
\mathcal{C}^{\otimes}_{\mathbf{x}}\to
\mathcal{D}^{\otimes}_{\mathbf{x}}$
is a lax $\mathcal{O}_i$-monoidal functor.
This implies that there is a commutative diagram
\[ \begin{array}{ccc}
    \mathcal{C}^{\otimes}_{(\mathbf{x},x)}
    &\stackrel{f_{(\mathbf{x},x)}}{\longrightarrow}&
    \mathcal{D}^{\otimes}_{(\mathbf{x},x)}\\[2mm]
    \bigg\downarrow&&\bigg\downarrow\\[2mm]
    \mathcal{C}^{\otimes}_{(\mathbf{x},x')}
    &\stackrel{f_{(\mathbf{x},x')}}{\longrightarrow}&
    \mathcal{D}^{\otimes}_{(\mathbf{x},x')}.\\[2mm]
   \end{array}\]
Since the vertical arrows are projections,
the above commutative diagram is left adjointable.
Hence we obtain a commutative diagram 
\[ \begin{array}{ccc}
    (\mathcal{C}^{\otimes})^{\vee}_{(\mathbf{x},x)}
    &\stackrel{g_{(\mathbf{x},x)}}{\longleftarrow}&
    (\mathcal{D}^{\otimes})^{\vee}_{(\mathbf{x},x)}\\[2mm]
    \bigg\downarrow&&\bigg\downarrow\\[2mm]
    (\mathcal{C}^{\otimes})^{\vee}_{(\mathbf{x},x')}
    &\stackrel{g_{(\mathbf{x},x')}}{\longleftarrow}&
    (\mathcal{D}^{\otimes})^{\vee}_{(\mathbf{x},x')}.\\[2mm]
   \end{array}\]
This means that $g_{\mathbf{x}}$ is an oplax 
$\mathcal{O}_i$-monoidal
functor. 
\qed

%\newpage
%%%\input{examples-o-p-duoidal}
%\newpage
%%%%\input{delta-product-duoidal}
%\newpage
%%%\input{lax-lax-formulation}
%\newpage
%%%\input{higher-monoidal}
%\newpage
%%%\input{o-product-duoidal-categories}
%\newpage
%%%\input{enriched_duoidal}

%\newpage

%{\footnotesize
%\input{ref}

%}

\end{document}